\documentclass[12pt]{amsart}

\usepackage{amssymb}
\usepackage[bbgreekl]{mathbbol}
\usepackage{graphicx}
\usepackage{hyperref}
\usepackage{ifthen}
\usepackage{mathrsfs}
\usepackage{pict2e}
\usepackage{xargs}
\usepackage{xspace}

\DeclareSymbolFontAlphabet{\amsmathbb}{AMSb}

\newcommand{\definedsymbol}[1]{$#1$}
\newcommand{\definedterm}[1]{\emph{#1}}

\newcommand{\Bairespace}[1][]{
  \ifthenelse{\equal{#1}{}}{\functions{\N}{\N}}{\functions{#1}{\N}}
}
\newcommand{\Bairetree}[1][]{
  \ifthenelse{\equal{#1}{}}{\functions{<\N}{\N}}{\functions{#1}{\N}}
}
\newcommand{\bbE}{\amsmathbb{E}}
\newcommand{\bbF}{\amsmathbb{F}}
\newcommand{\bbs}[1][]{\mathbb{s}_{#1}}
\newcommand{\calA}{\mathscr{A}}
\newcommand{\calI}{\mathcal{I}}
\newcommand{\calJ}{\mathcal{J}}
\newcommand{\calN}{\mathcal{N}}
\newcommand{\Cantorspace}[1][]{
  \ifthenelse{\equal{#1}{}}{\functions{\N}{2}}{\functions{#1}{2}}
}
\newcommand{\Cantortree}[1][]{
  \ifthenelse{\equal{#1}{}}{\functions{<\N}{2}}{\functions{#1}{2}}
}
\newcommand{\cardinality}[1]{|#1|}
\newcommand{\completeer}[1]{I(#1)}
\newcommand{\composition}{\circ}
\newcommandx{\concatenation}[2][1 = undefined, 2 = undefined]{
  \ifthenelse{\equal{#1}{undefined}}{{}\smallfrown}{
    \ifthenelse{\equal{#2}{undefined}}{\smallfrown_{#1}}{\bigoplus_{#1} #2}
  }
}
\newcommand{\constantsequence}[2]{\sequence{#1}^{#2}}
\newcommand{\DC}{\mathtt{DC}}
\newcommandx{\Deltaclass}[2][1=,2=]{
  \ifthenelse
    {\equal{#2}{}}
    {\mathbf{\Delta}_{#1}}
    {\mathbf{\Delta}^{#1}_{#2}}
}
\newcommand{\diagonal}[1]{\Delta(#1)}
\newcommandx{\disjointunion}[2][1 =, 2 =]{
  \ifthenelse{\equal{#1}{}}{\sqcup}{
    \ifthenelse{\equal{#2}{}}{\bigsqcup #1}{\bigsqcup_{#1} #2}
  }
}
\newcommand{\divides}{\mathrel{\mid}}
\newcommand{\domain}[1]{\mathrm{dom}(#1)}
\newcommand{\emptysequence}{\emptyset}
\newcommand{\equivalenceclass}[2]{[#1]_{#2}}
\newcommand{\extendedby}{\sqsubseteq}

\newcommand{\extension}[2][]{
  \ifthenelse{\equal{#1}{}}{\overline{#2}}{\overline{#2}_{#1}}
}
\newcommand{\extensions}[1]{\calN_{#1}}
\newcommand{\Ezero}{\bbE_0}
\newcommand{\from}{\colon}
\newcommandx{\functions}[3][3 =]{
  \ifthenelse{\equal{#3}{}}{#2^{#1}}{#2^{#1}_{#3}}
}
\newcommand{\Fzero}[1]{\bbF_{#1}}

\newcommand{\graph}[1]{\mathrm{graph}(#1)}
\newcommand{\Gzero}{\amsmathbb{G}_0}
\newcommand{\heightcorrection}[1]{\raisebox{0pt}[0pt][0pt]{#1}}
\newcommand{\horizontalsection}[2]{#1^{#2}}
\newcommand{\id}[1][]{\mathrm{id}_{#1}}
\newcommand{\image}[2]{#1(#2)}
\newcommandx{\intersection}[2][1 =, 2 =]{
  \ifthenelse{\equal{#1}{}}{\cap}{
    \ifthenelse{\equal{#2}{}}{\bigcap #1}{{\bigcap_{#1} #2}}
  }
}
\newcommandx{\interval}[3][3 =]{[#1, #2]_{#3}}
\newcommand{\into}{\hookrightarrow}
\newcommand{\inverse}[1]{#1^{-1}}
\newcommand{\length}[1]{|#1|}
\newcommand{\lifting}[1]{\tilde{#1}}
\newcommand{\mathand}{\text{ and }}

\renewcommand{\mod}[1]{\thinspace (\text{mod } #1)}
\newcommand{\N}{\amsmathbb{N}}

\newcommand{\onto}{\twoheadrightarrow}
\newcommand{\orbit}[2]{[#1]_{#2}}
\newcommand{\orbitequivalencerelation}[2]{E_{#1}^{#2}}
\newcommand{\pair}[2]{(#1, #2)}
\newcommand{\partialto}{\rightharpoonup}
\newcommandx{\Piclass}[2][1=,2=]{
  \ifthenelse{\equal{#2}{}}{\mathbf{\Pi}_{#1}}{\mathbf{\Pi}^{#1}_{#2}}
}
\newcommand{\positiveintegers}{\Z^+}
\newcommand{\preimage}[2]{#1^{-1}(#2)}
\newcommandx{\product}[2][1 =, 2 =]{
  \ifthenelse{\equal{#1}{}}{\times}{
    \ifthenelse{\equal{#2}{}}{\prod #1}{{\prod_{#1} #2}}
  }
}
\newcommand{\projection}[1][]{
  \ifthenelse{\equal{#1}{}}{\mathrm{proj}}{\mathrm{proj}_{#1}}
}
\newcommand{\quadruple}[4]{(#1, #2, #3, #4)}

\renewcommand{\restriction}[2]{#1 \upharpoonright #2}
\newcommand{\saturation}[2]{[#1]_{#2}}
\newcommandx{\sequence}[2][2 = undefined]{
  \ifthenelse{\equal{#2}{undefined}}{(#1)}{
    (#1)_{#2}
  }
}
\newcommandx{\set}[2][2 = undefined]{
  \ifthenelse{\equal{#2}{undefined}}{\{ #1 \}}{
    \{ #1 \suchthat #2 \}
  }
}
\newcommand{\setcomplement}[1]{\twiddle #1}
\newcommandx{\sets}[4][3 = undefined, 4 = undefined]{
  \ifthenelse{\equal{#3}{undefined}}{[#2]^{#1}}{
    \ifthenelse{\equal{#4}{undefined}}{[#2]^{#1}_{#3}}{[#2]^{#1}_{#3, #4}}
  }
}
\newcommandx{\Sigmaclass}[2][1=,2=]{
  \ifthenelse{\equal{#2}{}}
    {\mathbf{\Sigma}_{#1}}
    {\mathbf{\Sigma}^{#1}_{#2}}
}
\newcommand{\singleton}[1]{\set{#1}}

\newcommand{\strictlyextendedby}{\sqsubset}
\newcommand{\suchthat}{\mid}
\newcommand{\support}[1]{\mathrm{supp}(#1)}
\newcommand{\textexponent}[2]{$#1^{\text{#2}}$}
\newcommand{\triple}[3]{(#1, #2, #3)}
\newcommand{\twiddle}
  {\raisebox{1.5pt}{\scalebox{.75}{$\mathord{\sim}$}}}
\newcommandx{\union}[2][1 =, 2 =]{
  \ifthenelse{\equal{#1}{}}{\cup}{
    \ifthenelse{\equal{#2}{}}{\bigcup #1}{{\bigcup_{#1} #2}}
  }
}
\newcommand{\verticalsection}[2]{#1_{#2}}
\newcommand{\Z}{\amsmathbb{Z}}

\newcommand{\Baire}{Baire\xspace}
\newcommand{\Borel}{Bor\-el\xspace}
\newcommand{\Effros}{Eff\-ros\xspace}
\newcommand{\Feldman}{Feld\-man\xspace}
\newcommand{\Glimm}{Glimm\xspace}
\newcommand{\Hausdorff}{Haus\-dorff\xspace}
\newcommand{\Kechris}{Kech\-ris\xspace}
\newcommand{\Kuratowski}{Kur\-at\-ow\-ski\xspace}
\newcommand{\Lusin}{Lus\-in\xspace}
\newcommand{\Miller}{Mill\-er\xspace}
\newcommand{\Moore}{Moore\xspace}
\newcommand{\Novikov}{No\-vik\-ov\xspace}
\newcommand{\Polish}{Po\-lish\xspace}
\newcommand{\Solecki}{Sol\-eck\-i\xspace}
\newcommand{\Souslin}{Sous\-lin\xspace}
\newcommand{\Todorcevic}{To\-dor\-cev\-ic\xspace}
\newcommand{\Ulam}{U\-lam\xspace}

\newenvironment{lemmaproof}{
  
  \begin{proof}
}{\end{proof}}

\newenvironment{propositionproof}{
  
  \begin{proof}
}{\end{proof}}

\newenvironment{theoremproof}{
  
  \begin{proof}
}{\end{proof}}

\newtheorem{lemma}{Lemma}[section]
\newtheorem{proposition}[lemma]{Proposition}
\newtheorem{theorem}[lemma]{Theorem}

\newtheorem{introtheorem}{Theorem}

\theoremstyle{definition}
\newtheorem{remark}[lemma]{Remark}

\newtheorem{introremark}[introtheorem]{Remark}

\newtheorem*{acknowledgements}{Acknowledgements}

\begin{document}

%% Front matter

\begin{abstract}
  We establish generalizations of the \Feldman--\Moore theorem, the
  \Glimm--\Effros dichotomy, and the \Lusin--\Novikov uniformization
  theorem from \Polish spaces to their quotients by \Borel orbit equivalence
  relations.
\end{abstract}

\author{N. de Rancourt}

\address{
  N. de Rancourt \\
  Faculty of Mathematics and Physics \\
  Department of Mathematical Analysis \\
  Sokolovsk\'a 49/83 \\
  186 75 Praha 8 \\
  Czech Republic
 }

\email{rancourt@karlin.mff.cuni.cz}

\urladdr{
  \url{https://www2.karlin.mff.cuni.cz/~rancourt/}
}

\author{B.D. Miller}

\address{
  B.D. Miller \\
  Faculty of Mathematics \\
  University of Vienna \\
  Oskar Morgenstern Platz 1 \\
  1090 Wien \\
  Austria
 }

\email{benjamin.miller@univie.ac.at}

\urladdr{
  \url{https://homepage.univie.ac.at/benjamin.miller/}
}

\thanks{The authors were supported, in part, by FWF Grant P29999.}



\title[Feldman--Moore, Glimm--Effros, and Lusin--Novikov]
  {The Feldman--Moore, Glimm--Effros, and
    Lusin--Novikov theorems over quotients}

\maketitle

\section*{Introduction}

A topological space is \definedterm{\Polish} if it is separable and admits a
compatible complete metric. A subset of a topological space is \definedterm
{\Borel} if it is in the smallest $\sigma$-algebra containing the open sets.
Given an equivalence relation $E$ on a topological space $X$, we say that
a set $B \subseteq X / E$ is \definedterm{\Borel} if $\union[B]$ is \Borel.
More generally, given equivalence relations $E_i$ on topological spaces
$X_i$, we say that a set $R \subseteq \product[i \in I][X_i / E_i]$ is
\definedterm{weakly \Borel} if the corresponding set \heightcorrection
{$\lifting{R} = \set{\sequence{x_i}[i \in I] \in \product[i \in I][X_i]}[{\sequence
{\equivalenceclass{x_i}{E_i}}[i \in I] \in R}]$} is \Borel. Given equivalence
relations $E$ and $F$ on topological spaces $X$ and $Y$, we say that a
partial function $\phi \from X / E \partialto Y / F$ is \definedterm{strongly
\Borel} if its graph is weakly \Borel.

A \definedterm{$\sigma$-ideal} on a set $X$ is a family of subsets of $X$
that is closed under containment and countable unions. When $X$ is a
\Polish space, we say that an assignment $x \mapsto \calI_x$, sending
each point of $X$ to a $\sigma$-ideal on $X$, is \definedterm{strongly
\Borel-on-\Borel} if $\set{\pair{x}{y} \in X \times Y}[\verticalsection{R}{\pair
{x}{y}} \in \calI_x]$ is \Borel for all \Polish spaces $Y$ and \Borel sets $R
\subseteq (X \times Y) \times X$. A \Borel equivalence relation $E$ on a
\Polish space $X$ is \definedterm{strongly idealistic} if there is an
$E$-invariant strongly \Borel-on-\Borel assignment $x \mapsto \calI_x$
sending each point in $X$ to a $\sigma$-ideal on $X$ for which
$\equivalenceclass{x}{E} \notin \calI_x$. Following the usual abuse of
language, we say that an equivalence relation is \definedterm{countable} if
each of its equivalence classes is countable. The \Feldman--\Moore
theorem ensures that every countable \Borel equivalence relation on a
\Polish space is the orbit equivalence relation induced by a \Borel action of
a countable discrete group (see \cite[Theorem 1]{FeldmanMoore}), and the
proof of \cite[\S1.II.i]{Kechris:LocallyCompact} shows that every \Borel orbit
equivalence relation induced by a \Borel action of a \Polish group on a
\Polish space is strongly idealistic.

Our goal here is to generalize several basic results underlying the study
of countable \Borel equivalence relations from \Polish spaces to their
quotients by strongly idealistic \Borel equivalence relations. More precisely,
we provide countably-infinite bases of minimal counterexamples to such
generalizations.

A \definedterm{transversal} of an equivalence relation $E$ on $X$ is a set
$Y \subseteq X$ that intersects every $E$-class in exactly one point. An
\definedterm{embedding} of $E$ into an equivalence relation $F$ on a set
$Y$ is an injection $\pi \from X \to Y$ such that $w \mathrel{E} x \iff \pi(w)
\mathrel{F} \pi(x)$ for all $w, x \in X$. The \definedterm{diagonal} on $X$ is
given by $\diagonal{X} = \set{\pair{x}{y} \in X \times X}[x = y]$ and $\Ezero$
is the equivalence relation on $\Cantorspace$ given by $c \mathrel{\Ezero}
d \iff \exists n \in \N \forall m \ge n \ c(m) = d(m)$. We use $\Fzero{k}$ to
denote the index $k$ subequivalence relation of $\Ezero$ given by $c
\mathrel{\Fzero{k}} d \iff \exists n \in \N \forall m \ge n \ \sum_{i < m} c(i)
\equiv \sum_{i < m} d(i) \mod{k}$ for all $k \ge 2$. The following fact
generalizes the \Glimm--\Effros dichotomy for countable \Borel equivalence
relations from \Polish spaces to their quotients by strongly idealistic \Borel
equivalence relations:

\begin{introtheorem} \label{intro:transversal}
  Suppose that $X$ is the quotient of a \Polish space by a strongly idealistic
  \Borel equivalence relation and $E$ is a countable weakly \Borel
  equivalence relation on $X$. Then exactly one of the following holds:
  \begin{enumerate}
    \item The set $X$ is a countable union of \Borel transversals of $E$.
    \item There is a strongly \Borel embedding of $\Ezero / \bbF$ into $E$
      for some $\bbF \in \set{\diagonal{\Cantorspace}} \union \set{\Fzero{p}}
      [p \text{ is prime}]$.
  \end{enumerate}
\end{introtheorem}

The \definedterm{horizontal sections} of a set $R \subseteq X \times
Y$ are the sets of the form $\horizontalsection{R}{y} = \set{x \in X}[x
\mathrel{R} y]$ for $y \in Y$, whereas the \definedterm{vertical
sections} are those of the form $\verticalsection{R}{x} = \set{y \in Y}
[x \mathrel{R} y]$ for $x \in X$. The \definedterm{product} of equivalence
relations $E$ on $X$ and $F$ on $Y$ is the equivalence relation on $X
\times Y$ given by $\pair{x}{y} \mathrel{E \times F} \pair{x'}{y'} \iff (x
\mathrel{E} x' \mathand y \mathrel{F} y')$. A \definedterm{rectangular
homomorphism} from a binary relation $R$ on $X \times Y$ to a binary
relation $R'$ on $X' \times Y'$ is a function of the form $\phi \times \psi$,
where $\phi \from X \to X'$ and $\psi \from Y \to Y'$, with the property that
$\image{(\phi \times \psi)}{R} \subseteq R'$. The following fact generalizes
the \Lusin--\Novikov uniformization theorem from \Polish spaces to their
quotients by strongly idealistic \Borel equivalence relations:

\begin{introtheorem} \label{intro:uniformization}
  Suppose that $X$ is the quotient of a \Polish space by a \Borel
  equivalence relation, $Y$ is the quotient of a \Polish space by a
  strongly idealistic \Borel equivalence relation, and $R \subseteq X \times
  Y$ is a weakly \Borel set whose vertical sections are countable. Then
  exactly one of the following holds:
  \begin{enumerate}
    \item The set $\image{\projection[X]}{R}$ is \Borel and there are strongly
      \Borel functions $\phi_n \from \image{\projection[X]}{R} \to Y$ for which
      $R = \union[n \in \N][\graph{\phi_n}]$.
    \item There is an injective strongly \Borel rectangular homomorphism
      from $\Ezero / (\Ezero \times \bbF)$ to $R$ for some $\bbF \in \set
      {\diagonal{\Cantorspace}} \union \set{\Fzero{p}}[p \text{ is prime}]$.
  \end{enumerate}
\end{introtheorem}

\begin{introremark} \label{intro:uniformization:E}
  In the special case that $X$ is a \Polish space, condition (2) cannot hold,
  so condition (1) always holds.
\end{introremark}

\begin{introremark} \label{intro:uniformization:F}
  In the special case that $Y$ is a \Polish space, condition (2) cannot hold
  when $\bbF \in \set{\Fzero{p}}[p \text{ is prime}]$, so $\Ezero / (\Ezero
  \times \diagonal{\Cantorspace})$ is the minimal counterexample to
  condition (1), answering a question that originally arose in a conversation
  with \Kechris.
\end{introremark}

The \definedterm{complete equivalence relation} on $X$ is given by
$\completeer{X} = X \times X$, and the \definedterm{disjoint union} of
equivalence relations $E_0$ and $E_1$ on $X$ is the equivalence relation
on $X \times 2$ given by $\pair{x}{i} \mathrel{E_0 \disjointunion E_1} \pair
{y}{j} \iff (i = j \mathand x \mathrel{E_i} y)$. The following fact generalizes
the \Feldman--\Moore theorem from \Polish spaces to their quotients by
strongly idealistic \Borel equivalence relations:

\begin{introtheorem} \label{intro:graph}
  Suppose that $X$ is the quotient of a \Polish space by a strongly idealistic
  \Borel equivalence relation and $E$ is a countable weakly \Borel
  equivalence relation on $X$. Then exactly one of the following holds:
  \begin{enumerate}
    \item There is a countable group $\Gamma$ of strongly \Borel
      automorphisms of $X$ for which $E = \orbitequivalencerelation
      {\Gamma}{X}$.
    \item There is a strongly \Borel embedding of $(\Ezero \times \completeer
      {2}) / (\Ezero \disjointunion \bbF)$ into $E$ for some $\bbF \in \set
      {\diagonal{\Cantorspace}} \union \set{\Fzero{p}}[p \text{ is prime}]$.
  \end{enumerate}
\end{introtheorem}

In \S\ref{gzero}, we provide the original classical proof of the $\Gzero$
dichotomy. While the underlying argument has been known for some
time, it has somehow never appeared in full. We provide it here for
future reference, and because our later arguments depend not only
upon the $\Gzero$ dichotomy, but its proof, as well.

In \S\ref{generalization}, we establish the special case of Theorem
\ref{intro:uniformization} referred to in Remark \ref
{intro:uniformization:E}. While it is not strictly necessary to establish
this case separately, our argument is substantially simpler than those
appearing in the latter sections---to which it serves as a stepping
stone---and also yields a more general result.

In \S\ref{finitebasis}, we establish versions of the special cases of Theorem
\ref{intro:transversal} where every $E$-class has cardinality at most $k$,
Theorem \ref{intro:uniformization} where every vertical section of $R$ 
has cardinality at most $k$, and Theorem \ref{intro:graph} where every
$E$-class has cardinality at most $k + 1$. In these special cases, condition
(2) can only hold when $\bbF \in \set{\Fzero{p}}[p \le k]$.

In \S\ref{dichotomy}, we establish analogs of Theorems \ref
{intro:transversal}, \ref{intro:uniformization}, and \ref{intro:graph} in which
the first conditions are relaxed to allow for bounded finite error, but the
second conditions are strengthened to ensure that $\bbF = \diagonal
{\Cantorspace}$. While it is not difficult to derive these from \cite[Theorem
1]{CCM} and the results of \S\ref{generalization}, we instead show them
using a substantial simplification of the proof of the former. The special
case of Theorem \ref{intro:uniformization} referred to in Remark \ref
{intro:uniformization:F} easily follows, and by combining the results of the
final three sections, we obtain Theorems \ref{intro:transversal}, \ref
{intro:uniformization}, and \ref{intro:graph}.

\section{The $\Gzero$ dichotomy} \label{gzero}

Endow $\N$ with the discrete topology and $\Bairespace$ with the
corresponding product topology. A topological space is \definedterm
{analytic} if it is a continuous image of a closed subset of $\Bairespace$,
and a subset of a topological space is \definedterm{co-analytic} if its
complement is analytic. Every \Polish space is analytic (see, for example,
\cite[Theorem 7.9]{Kechris}), and \Souslin's theorem ensures that a subset
of an analytic \Hausdorff space is \Borel if and only if it is analytic and
co-analytic (see, for example, \cite[Theorem 14.11]{Kechris}\footnote{While
the results of \cite{Kechris} are stated for \Polish spaces, the proofs of
those to which we refer go through just as easily in the generality discussed
here.}).

Given a natural number $n \ge 2$, a sequence $\sequence{X_i}[i < n]$,
and sets $R \subseteq \product[i < n][X_i]$ and $Y_i \subseteq X_i$ for all
$i < n$, we say that $\sequence{Y_i}[i < n]$ is \definedterm
{$R$-independent} if $R \intersection \product[i < n][Y_i] = \emptyset$. We
identify $\product[i < n][X_i]$ with $\product[i \le m][X_i] \times \product[m
< i < n][X_i]$ for all $m < n$.

\begin{proposition} \label{gzero:separation:relation}
  Suppose that $n \ge 2$ is a natural number, $\sequence{X_i}[i < n]$
  is a sequence of \Hausdorff spaces, $A_i \subseteq X_i$ is analytic
  for all $i < n$, $R \subseteq \product[i < n][X_i]$ is analytic, and
  $\sequence{A_i}[i < n]$ is $R$-independent. Then there are \Borel
  sets $B_i \subseteq X_i$ such that $A_i \subseteq B_i$ for all $i < n$
  and $\sequence{B_i}[i < n]$ is $R$-independent.
\end{proposition}

\begin{propositionproof}
  It is sufficient to show that if $m < n$, $B_i \subseteq X_i$ is a \Borel
  set such that $A_i \subseteq B_i$ for all $i < m$, and
  $\sequence{B_i}[i < m] \union \sequence{A_i}[m \le i < n]$ is
  $R$-independent, then there is a \Borel set $B_m \subseteq X_m$
  for which $A_m \subseteq B_m$ and $\sequence{B_i}[i
  \le m] \union \sequence{A_i}[m < i < n]$ is $R$-independent.
  Towards this end, let \definedsymbol{\projection[i]} denote the
  projection function from $\functions{n}{X}$ onto the \textexponent{i}
  {th} coordinate for all $i < n$. As the class of analytic \Hausdorff
  spaces is closed under continuous images, intersections, products,
  and \Borel subsets (see, for example, \cite[Proposition
  14.4]{Kechris}), the fact that $R \intersection ((\product[i < m][B_i])
  \times X \times (\product[m < i < n][A_i])) = R \intersection ((\product
  [i \le m][\image{\projection[i]}{R}]) \times (\product[m < i < n][A_i]))
  \intersection ((\product[i < m][B_i]) \times \functions{n - m}{X})$
  ensures that $\image{\projection[m]}{R \intersection ((\product[i < m]
  [B_i]) \times X \times (\product[m < i < n][A_i]))}$ is analytic, and
  since $A_m$ is disjoint from this set, \Lusin's separation theorem
  (see, for example, \cite[Theorem 14.7]{Kechris}) yields a \Borel set
  $B_m \subseteq X$ containing $A_m$ and disjoint from $\image
  {\projection[m]}{R \intersection ((\product[i < m][B_i]) \times X \times
  (\product[m < i < n][A_i]))}$, thus $\sequence{B_i}[i \le m] \union
  \sequence{A_i}[m < i < n]$ is $R$-independent.
\end{propositionproof}

For all sets $X$ and natural numbers $n$, we use \definedsymbol
{\constantsequence{X}{n}} to denote the constant sequence of length $n$
with value $X$. For all $n \ge 2$, an \definedterm{$n$-dimensional
dihypergraph} on a set $X$ is an $n$-ary relation $G$ on $X$ consisting
solely of non-constant sequences. We say that a set $Y \subseteq X$ is
\definedterm{$G$-independent} if $\constantsequence{Y}{n}$ is
$G$-independent.

\begin{proposition} \label{gzero:separation:graph}
  Suppose that $n \ge 2$, $X$ is a \Hausdorff
  space, $G$ is an analytic $n$-dimensional dihypergraph on $X$,
  and $A \subseteq X$ is a $G$-independent analytic set. Then there
  is a $G$-independent \Borel set $B \subseteq X$ for which $A
  \subseteq B$.
\end{proposition}

\begin{propositionproof}
  By Proposition \ref{gzero:separation:relation}, there are \Borel sets
  $B_i \subseteq X$ such that $A \subseteq B_i$ for all $i < n$ and
  $\sequence{B_i}[i < n]$ is $G$-independent. Set $B = \intersection
  [i < n][B_i]$.
\end{propositionproof}

We use $\functions{<\N}{X}$ to denote $\union[n \in \N][\functions{n}
{X}]$, \definedsymbol{\sequence{i}} to denote the singleton sequence
with value $i$, \definedsymbol{\extendedby} to denote extension, and
\definedsymbol{\negmedspace \concatenation} to denote concatenation of
sequences. Following standard practice, we also use $\extensions{s}$ to
denote $\set{b \in \Bairespace}[s \extendedby b]$ or $\set{c \in
\Cantorspace}[s \extendedby c]$ (with the context determining which of the
two we have in mind). A \definedterm{digraph} is a two-dimensional
dihypergraph. Fix sequences $\bbs[n] \in \Cantorspace[n]$ such that
$\forall s \in \Cantortree \exists n \in \N \ s \extendedby \bbs[n]$ and
let $\Gzero$ denote the digraph on $\Cantorspace$ given by $\Gzero =
\set{\sequence{\bbs[n] \concatenation \sequence{i} \concatenation c}
[i < 2]}[c \in \Cantorspace \mathand n \in \N]$.

\begin{proposition} \label{gzero:meager}
  Every $\Gzero$-independent set $B \subseteq \Cantorspace$ with
  the \Baire property is meager.
\end{proposition}

\begin{propositionproof}
  Suppose, towards a contradiction, that $B$ is not meager, and fix a
  sequence $s \in \Cantortree$ for which $B$ is comeager in
  $\extensions{s}$ (see, for example, \cite[Proposition 8.26]{Kechris})
  and $n \in \N$ for which $s \extendedby \bbs[n]$, and let $\iota$ be
  the involution of $\extensions{\bbs[n]}$ given by $\iota(\bbs[n]
  \concatenation \sequence{0} \concatenation c) = \bbs[n]
  \concatenation \sequence{1} \concatenation c$ for all $c \in
  \Cantorspace$. As $\iota$ is a homeomorphism, it follows that $B
  \intersection \image{\iota}{B}$ is comeager in $\extensions{\bbs[n]}$
  (see, for example, \cite[Exercise 8.45]{Kechris}), so $B \intersection
  \image{\iota}{B} \intersection \extensions{\bbs[n] \concatenation
  \sequence{0}} \neq \emptyset$. But $\pair{c}{\iota(c)} \in \Gzero
  \intersection (B \times B)$ for all $c \in B \intersection \image{\iota}
  {B} \intersection \extensions{\bbs[n] \concatenation \sequence{0}}$,
  contradicting the fact that $B$ is $\Gzero$-independent.
\end{propositionproof}

An \definedterm{$I$-coloring} of a digraph $G$ on a set $X$ is a function
$c \from X \to I$ such that $c(x) \neq c(y)$ for all $\pair{x}{y} \in G$, or
equivalently, such that $\preimage{c}{\set{i}}$ is $G$-independent for all $i
\in I$. A \definedterm{homomorphism} from a binary relation $R$ on a set
$X$ to a binary relation $S$ on a set $Y$ is a function $\phi \from X \to Y$
for which $\image{(\phi \times \phi)}{R} \subseteq S$. 

\begin{theorem}[\Kechris--\Solecki--\Todorcevic]
  \label{gzero:main}
  Suppose that $X$ is a \Hausdorff space and $G$ is an analytic
  digraph on $X$. Then exactly one of the following holds:
  \begin{enumerate}
    \setlength{\itemindent}{-4.35pt}
    \item There is a \Borel $\N$-coloring $c \from X \to \N$ of $G$.
    \item There is a continuous homomorphism $\phi \from
      \Cantorspace \to X$ from $\Gzero$ to $G$.
  \end{enumerate}
\end{theorem}

\begin{theoremproof}[Proof (\Miller)]
  To see that the two conditions are mutually exclusive, note that if
  both hold, then $c \composition \phi$ is a \Borel $\N$-coloring of
  $\Gzero$, so there exists $n \in \N$ for which $\preimage{(c
  \composition \phi)}{\set{n}}$ is a non-meager $\Gzero$-independent
  \Borel set, contradicting Proposition \ref{gzero:meager}.

  It remains to show that at least one of the two conditions holds. We
  can clearly assume that $G$ is not empty, so there are continuous
  surjections $\phi_G \from \Bairespace \onto G$ and $\phi_X \from
  \Bairespace \onto \union[i < 2][\image{\projection[i]}{G}]$.
   
  We will recursively construct a decreasing sequence $\sequence
  {B^\alpha}[\alpha < \omega_1]$ of \Borel subsets of $X$, off of which
  there are \Borel $\N$-colorings of $G$. We begin by setting $B^0 =
  X$ and we define $B^\lambda = \intersection[\alpha < \lambda][B^\alpha]$
  for all limit ordinals $\lambda < \omega_1$. To describe the construction
  of $B^{\alpha + 1}$ from $B^\alpha$, we require several preliminaries.
  
  An \definedterm{approximation} is a triple of the form $a = \triple{n^a}
  {\phi^a}{\sequence{\psi_n^a}[n < n^a]}$, where $n^a \in \N$, $\phi^a
  \from \Cantorspace[n^a] \to \Bairespace[n^a]$, and $\psi_n^a \from
  \Cantorspace[n^a - (n+1)] \to \Bairespace[n^a]$ for all $n < n^a$. A
  \definedterm{one-step extension} of $a$ is an approximation $b$
  such that:
  \begin{enumerate}
    \renewcommand{\theenumi}{\alph{enumi}}
    \setlength{\itemindent}{-5.7pt}
    \item $n^b = n^a + 1$,
    \item $\forall s \in \Cantorspace[n^a] \forall t \in \Cantorspace[n^b]
      \ (s \extendedby t \implies \phi^a(s) \extendedby
        \phi^b(t))$, and
    \item $\forall n < n^a \forall s \in \Cantorspace[n^a - (n + 1)] \forall t
      \in \Cantorspace[n^b - (n + 1)] \ (s \extendedby t \implies
        \psi^a_n(s) \extendedby \psi^b_n(t))$.
  \end{enumerate}
  A \definedterm{configuration} is a triple of the form $\gamma = \triple
  {n^\gamma}{\phi^\gamma}{\sequence{\psi_n^\gamma}[n <
  n^\gamma]}$, where $n^\gamma \in \N$, $\phi^\gamma \from
  \Cantorspace[n^\gamma] \to \Bairespace$, $\psi_n^\gamma \from
  \Cantorspace[n^\gamma - (n+1)] \to \Bairespace$ for all $n <
  n^\gamma$, and $(\phi_G \composition \psi_n^\gamma)(t) = \pair
  {(\phi_X \composition \phi^\gamma)(\bbs[n] \concatenation \sequence
  {0} \concatenation t)}{(\phi_X \composition \phi^\gamma)(\bbs[n]
  \concatenation \sequence{1} \concatenation t)}$ for all $n <
  n^\gamma$ and $t \in \Cantorspace[n^\gamma - (n + 1)]$.
  We say that $\gamma$ is \definedterm{compatible} with a set $X'
  \subseteq X$ if $\image{(\phi_X \composition \phi^\gamma)}
  {\Cantorspace[n^\gamma]} \subseteq X'$, and \definedterm
  {compatible} with $a$ if:
  \begin{enumerate}
    \renewcommand{\theenumi}{\roman{enumi}}
    \item $n^a = n^\gamma$,
    \item $\forall t \in \Cantorspace[n^a] \ \phi^a(t) \extendedby
      \phi^\gamma(t)$, and
    \item $\forall n < n^a \forall t \in \Cantorspace[n^a - (n + 1)]
      \ \psi^a_n(t) \extendedby \psi^\gamma_n(t)$.
  \end{enumerate}
   An approximation $a$ is \definedterm{$X'$-terminal} if no
   configuration is compatible with both $X'$ and a one-step extension
   of $a$. Let $A(a, X')$ denote the set of points of the form $(\phi_X
   \composition \phi^\gamma)(\bbs[n^a])$, where $\gamma$ varies
   over all configurations compatible with $a$ and $X'$.
  
  \begin{lemma} \label{gzero:terminalimpliesindependent}
    Suppose that $X' \subseteq X$ and $a$ is an $X'$-terminal
    approximation. Then $A(a, X')$ is $G$-independent.
  \end{lemma}
  
  \begin{lemmaproof}
    Suppose, towards a contradiction, that there are configurations
    $\gamma_0$ and $\gamma_1$, both compatible with $a$ and $X'$,
    with the property that $\pair{(\phi_X \composition \phi^{\gamma_0})
    (\bbs[n^a])}{(\phi_X \composition \phi^{\gamma_1})(\bbs[n^a])} \in
    G$. Fix a sequence $d \in \Bairespace$ such that $\phi_G(d) = \pair
    {(\phi_X \composition \phi^{\gamma_0})(\bbs[n^a])}{(\phi_X
    \composition \phi^{\gamma_1})(\bbs[n^a])}$, and let $\gamma$ be
    the configuration given by $n^\gamma = n^a + 1$, $\phi^\gamma(t
    \concatenation \sequence{i}) = \phi^{\gamma_i}(t)$ for all $i < 2$
    and $t \in \Cantorspace[n^a]$, $\psi_n^\gamma(t \concatenation
    \sequence{i}) = \psi_n^{\gamma_i}(t)$ for all $i < 2$, $n < n^a$,
    and $t \in \Cantorspace[n^a - (n + 1)]$, and $\psi_{n^a}^\gamma
    (\emptysequence) = d$. Then $\gamma$ is compatible with a
    one-step extension of $a$, contradicting the fact that $a$ is
    $X'$-terminal.
  \end{lemmaproof}
  
  For all $B^\alpha$-terminal approximations $a$, Proposition \ref
  {gzero:separation:graph} yields a $G$-independent \Borel set $B(a,
  B^\alpha) \supseteq A(a, B^\alpha)$. Let $B^{\alpha +1}$ be the set
  obtained from $B^\alpha$ by subtracting the union of the sets of the
  form $B(a, B^\alpha)$, where $a$ varies over all $B^\alpha$-terminal
  approximations.
  
  \begin{lemma}
    \label{gzero:nonterminalapproximationextension}
    Suppose that $\alpha < \omega_1$ and $a$ is a
    non-$B^{\alpha + 1}$-terminal approximation. Then $a$ has a
    non-$B^\alpha$-terminal one-step extension.
  \end{lemma}
  
  \begin{lemmaproof}
    Fix a one-step extension $b$ of $a$ for which there is a
    configuration $\gamma$ compatible with $b$ and $B^{\alpha + 1}$.
    Then $(\phi_X \composition \phi^\gamma)(\bbs[n^b]) \in B^{\alpha +
    1}$, so $b$ is not $B^\alpha$-terminal. 
  \end{lemmaproof}
  
  Fix $\alpha < \omega_1$ such that the families of $B^\alpha$- and
  $B^{\alpha + 1}$-terminal approximations coincide, and let $a_0$ be
  the unique approximation for which $n^{a_0} = 0$. As $A(a_0, X') =
  \union[i < 2][\image{\projection[i]}{G}] \intersection X'$ for all $X'
  \subseteq X$, we can assume that $a_0$ is not $B^\alpha$-terminal,
  since otherwise $B^{\alpha+1}$ is disjoint from $\union[i < 2][\image
  {\projection[i]}{G}]$, so there is a \Borel $\N$-coloring of $G$.
  
  By recursively applying Lemma \ref
  {gzero:nonterminalapproximationextension}, we obtain
  non-$B^\alpha$-terminal one-step extensions $a_{n+1}$ of $a_n$
  for all $n \in \N$. Define $\phi, \psi_n \from \Cantorspace \to
  \Bairespace$ by $\phi(c) = \union[n \in \N][\phi^{a_n}(\restriction{c}
  {n})]$ and $\psi_n(c) = \union[m > n][\psi_n^{a_m}(\restriction{c}{(m
  - (n + 1))})]$ for all $n \in \N$. Clearly these functions are continuous.
  
  To establish that the function $\pi = \phi_X \composition \phi$ is a
  homomorphism from $\Gzero$ to $G$, we will show the stronger
  fact that if $c \in \Cantorspace$ and $n \in \N$, then
  \begin{equation*}
    (\phi_G \composition \psi_n)(c) = \pair{(\phi_X \composition \phi)
      (\bbs[n] \concatenation \sequence{0} \concatenation c)}{(\phi_X
        \composition \phi)(\bbs[n] \concatenation \sequence{1}
          \concatenation c)}.
  \end{equation*}
  As $X \times X$ is \Hausdorff, it is sufficient to show that if $U$ is an
  open neighborhood of $\pair{(\phi_X \composition \phi)(\bbs[n]
  \concatenation \sequence{0} \concatenation c)}{(\phi_X \composition
  \phi)(\bbs[n] \concatenation \sequence{1} \concatenation c)}$ and
  $V$ is an open neighborhood of $(\phi_G \composition \psi_n)(c)$,
  then $U \intersection V \neq \emptyset$. Towards this end, fix $m >
  n$ such that $\image{\phi_X}{\extensions{\phi^{a_m}(\bbs[n]
  \concatenation \sequence{0} \concatenation s)}} \times \image
  {\phi_X}{\extensions{\phi^{a_m}(\bbs[n] \concatenation \sequence{1}
  \concatenation s)}} \subseteq U$ and $\image{\phi_G}{\extensions
  {\psi_n^{a_m}(s)}} \subseteq V$ where $s = \restriction{c}{(m - (n +
  1))}$. The fact that $a_m$ is not $B^\alpha$-terminal yields a
  configuration $\gamma$ compatible with $a_m$. Then $\pair{(\phi_X
  \composition \phi^\gamma)(\bbs[n] \concatenation \sequence{0}
  \concatenation s)}{(\phi_X \composition \phi^\gamma)(\bbs[n]
  \concatenation \sequence{1} \concatenation s)} \in U$ and $(\phi_G
  \composition \psi_n^\gamma)(s) \in V$, so $U \intersection V \neq
  \emptyset$.
\end{theoremproof}

\begin{remark} \label{gzero:DC}
  The apparent use of choice beyond $\DC$ in the above argument
  can be eliminated by first running the simplified version without
  Proposition \ref{gzero:separation:graph}, i.e., by setting $B(a,
  B^\alpha) = A(a, B^\alpha)$, in order to obtain an upper bound
  $\alpha' < \omega_1$ on the least ordinal $\alpha < \omega_1$ for
  which the sets of $B^\alpha$- and $B^{\alpha+1}$-terminal
  approximations coincide.
\end{remark}

\section{Generalizations} \label{generalization}

Given an equivalence relation $E$ on a set $X$, the \definedterm
{$E$-saturation} of a set $Y \subseteq X$ is given by $\saturation{Y}
{E} = \set{x \in X}[\exists y \in Y \ x \mathrel{E} y]$. We say that an
$n$-dimensional dihypergraph $G$ on $X$ is \definedterm
{$E$-invariant} if $x \in G \iff y \in G$ for all $x, y \in \functions
{n}{X}$ with the property that $\forall i < n \ x(i) \mathrel{E} y(i)$.

\begin{proposition} \label{generalization:separation}
  Suppose that $n \ge 2$, $X$ is a \Hausdorff
  space, $E$ is an analytic equivalence relation on $X$, $G$ is an
  $E$-invariant analytic $n$-dimensional dihypergraph on $X$, and $A
  \subseteq X$ is a $G$-independent analytic set. Then there is an
  $E$-invariant $G$-independent \Borel set $B \subseteq X$ for which
  $A \subseteq B$.
  \end{proposition}

\begin{propositionproof}
  Set $A_0 = A$. Given $n \in \N$ and a $G$-independent analytic set
  $A_n \subseteq X$, Proposition \ref{gzero:separation:graph} yields a
  $G$-independent \Borel set $B_n \subseteq X$ containing $A_n$,
  and since the class of analytic \Hausdorff spaces is closed under
  continuous images and \Borel subsets, it follows that the
  $G$-independent set $A_{n+1} = \saturation{B_n}{E} = \image
  {\projection[0]}{E \intersection (X \times B_n)}$ is analytic. Define $B
  = \union[n \in \N][A_n] = \union[n \in \N][B_n]$.
\end{propositionproof}

Clearly $\Gzero \subseteq \Ezero$. The following observation is
slightly less obvious:

\begin{proposition} \label{generalization:generates}
  The smallest equivalence relation $E$ on $\Cantorspace$ containing
  $\Gzero$ is $\Ezero$.
\end{proposition}

\begin{propositionproof}
  To see that $\forall c \in \Cantorspace \forall u, v \in \Cantorspace[n]
  \ u \concatenation c \mathrel{E} v \concatenation c$ for all $n \in \N$,
  observe that if it holds at some $n \in \N$, $c \in \Cantorspace$, and
  $u, v \in \Cantorspace[n]$, then $u \concatenation \sequence{0}
  \concatenation c \mathrel{E} \bbs[n] \concatenation \sequence{0}
  \concatenation c \mathrel{\Gzero} \bbs[n] \concatenation \sequence{1}
  \concatenation c \mathrel{E} v \concatenation \sequence{1}
  \concatenation c$, so it holds at $n + 1$.
\end{propositionproof}

We next note another connection between \Baire category and
$\Gzero$:

\begin{proposition} \label{generalization:meager}
  Suppose that $E$ and $F$ are equivalence relations on
  $\Cantorspace$ with the \Baire property, every $E$-class is a
  countable union of $(E \intersection F)$-classes, and $F \intersection
  \Gzero = \emptyset$. Then $E$ and $F$ are meager.
\end{proposition}

\begin{propositionproof}
  Suppose, towards a contradiction, that $F$ is not meager. As $F$
  has the \Baire property, the \Kuratowski--\Ulam theorem (see, for
  example, \cite[Theorem 8.41]{Kechris}) yields a sequence $c \in
  \Cantorspace$ whose $F$-class has the \Baire property but is not
  meager, in which case Proposition \ref{gzero:meager} yields a pair
  $\pair{a}{b} \in \restriction{\Gzero}{\equivalenceclass{c}{F}}$,
  contradicting the fact that $F \intersection \Gzero = \emptyset$. It
  therefore follows that $F$ is meager.

  The \Kuratowski--\Ulam theorem now ensures that every $F$-class
  is meager, in which case every $(E \intersection F)$-class is meager,
  so every $E$-class is meager, thus $E$ is meager.
\end{propositionproof}

A \definedterm{homomorphism} from a sequence $\sequence{R_i}[i \in I]$
of binary relations on a set $X$ to a sequence $\sequence{S_i}[i \in I]$ of
binary relations on a set $Y$ is a function $\phi \from X \to Y$ that is a
homomorphism from $R_i$ to $S_i$ for all $i \in I$. In this section, we only
need the weakening of the following fact where $\pair{\Fzero{k}}{\Ezero
\setminus \Fzero{k}}$ is replaced with $\Ezero$:

\begin{proposition} \label{generalization:homomorphism}
  Suppose that $k \ge 2$, $D$ is a closed nowhere-dense binary relation
  on $\Cantorspace$, and $R$ is a meager binary relation on
  $\Cantorspace$. Then there is a continuous homomorphism $\phi \from
  \Cantorspace \to \Cantorspace$ from $\quadruple{\setcomplement
  {\diagonal{\Cantorspace}}}{\Fzero{k}}{\Ezero \setminus \Fzero{k}}
  {\setcomplement{\Ezero}}$ to $\quadruple{\setcomplement{D}}{\Fzero{k}}
  {\Ezero \setminus \Fzero{k}}{\setcomplement{R}}$.
\end{proposition}

\begin{propositionproof}
  Fix a decreasing sequence $\sequence{U_n}[n \in \N]$ of dense open
  symmetric subsets of $\setcomplement{D}$ whose intersection is disjoint
  from $R$, as well as $\phi_0 \from \Cantorspace[0] \to \Cantorspace[0]$.
  
  \begin{lemma} \label{basis:F0:extension}
    Suppose that $n \in \N$ and $\phi_n \from \Cantorspace[n] \to
    \Cantortree$. Then there exist $\ell_n > 0$ and $u_n \in \Cantorspace
    [\ell_n] \times \Cantorspace[\ell_n]$ with the property that $1 + \sum_{j <
    \ell_n} u_n(0)(j) \equiv \sum_{j < \ell_n} u_n(1)(j) \mod{k}$ and $\product
    [i < 2][\extensions{\phi_{n+1}({t(i)} \concatenation \sequence{i})}]
    \subseteq U_n$ for all $t \in \Cantorspace[n] \times \Cantorspace[n]$,
    where $\phi_{n+1} \from \Cantorspace[n+1] \to \Cantortree$ is given by
    $\phi_{n+1}(t \concatenation \sequence{i}) = \phi_n(t) \concatenation
    u_n(i)$ for all $i < 2$ and $t \in \Cantorspace[n]$.
  \end{lemma}
  
  \begin{lemmaproof}
    Fix an enumeration $\sequence{t_j}[j < 4^n]$ of $\Cantorspace[n] \times
    \Cantorspace[n]$ and $u_{0,n} \in \Cantortree \times \Cantortree$. Given
    $j < 4^n$ and $u_{j,n} \in \Cantortree \times \Cantortree$, fix $u_{j+1,n}
    \in \Cantortree \times \Cantortree$ such that:
    \begin{enumerate}
      \item $\forall i < 2 \ u_{j,n}(i) \extendedby u_{j+1,n}(i)$.
      \item $\product[i < 2][\extensions{\phi_n(t_j(i)) \concatenation
        u_{j+1,n}(i)}] \subseteq U_n$.
    \end{enumerate}
    Then any $\ell_n > 0$ and $u_n \in \Cantorspace[\ell_n] \times
    \Cantorspace[\ell_n]$ such that $\forall i < 2 \ u_{4^n,n}(i) \extendedby
    u_n(i)$ and $1 + \sum_{j < \ell_n} u_n(0)(j) \equiv \sum_{j < \ell_n} u_n
    (1)(j) \mod{k}$ are as desired.
  \end{lemmaproof}
  
  As $\phi_n(t) \strictlyextendedby \phi_{n+1}(t \concatenation \sequence
  {i})$ for all $i < 2$, $n \in \N$, and $t \in \Cantorspace[n]$, we obtain a
  continuous function $\phi \from \Cantorspace \to \Cantorspace$ by setting
  $\phi(c) = \union[n \in \N][\phi_n(\restriction{c}{n})]$ for all $c \in
  \Cantorspace$. To see that $\phi$ is a homomorphism from
  $\setcomplement{\diagonal{\Cantorspace}}$ to $\setcomplement{D}$,
  note that if $c \in \setcomplement{\diagonal{\Cantorspace}}$, then there
  exists $n \in \N$ for which $c(0)(n) \neq c(1)(n)$, so $\sequence{\phi(c(i))}
  [i < 2] \in \product[i < 2][\extensions{\phi_{n+1}(\restriction{c(i)}{(n+1)})}]
  \subseteq U_n \subseteq \setcomplement{D}$. To see that $\phi$ is a
  homomorphism from $\pair{\Fzero{k}}{\Ezero \setminus \Fzero{k}}$ to
  $\pair{\Fzero{k}}{\Ezero \setminus \Fzero{k}}$, note that if $c \in \Ezero$,
  then there exists $n \in \N$ such that $c(0)(m) = c(1)(m)$ for all $m \ge
  n$, and setting $\ell = \sum_{j < n} \ell_j$, it follows that $\phi(c(0))(m) =
  \phi(c(1))(m)$ for all $m \ge \ell$, in which case
  \begin{align*}
    c&(0) \mathrel{\Fzero{k}} c(1) \\
      & \textstyle \iff \sum_{m < n} c(0)(m) \equiv \sum_{m < n} c(1)(m)
        \mod{k} \\
      & \textstyle \iff \sum_{m < n} \sum_{j < \ell_m} u_m(c(0))(j) \equiv
        \sum_{m < n} \sum_{j < \ell_m} u_m(c(1))(j) \mod{k} \\
      & \textstyle \iff \sum_{m < \ell} \phi_n(\restriction{c(0)}{n})(m) \equiv
        \sum_{m < \ell} \phi_n(\restriction{c(1)}{n})(m) \mod{k} \\
      & \iff \phi(c(0)) \mathrel{\Fzero{k}} \phi(c(1)).
  \end{align*}
  To see that $\phi$ is a homomorphism from $\setcomplement{\Ezero}$ to
  $\setcomplement{R}$, note that if $c \in \setcomplement{\Ezero}$, then
  there are infinitely many $n \in \N$ with the property that $\sequence{\phi
  (c(i))}[i < 2] \in \product[i < 2][\extensions{\phi_{n+1}(\restriction{c(i)}
  {(n+1)})}] \subseteq U_n$, so $\sequence{\phi(c(i))}[i < 2] \in
  \setcomplement{R}$.
\end{propositionproof}

A \definedterm{partial transversal} of an equivalence relation $E$ on a set
$X$ over a subequivalence relation $F$ is a set $Y \subseteq X$ for which
$\restriction{E}{Y} = \restriction{F}{Y}$.

\begin{theorem} \label{generalization:transversal}
  Suppose that $X$ is a \Hausdorff space, $E$ is an analytic equivalence
  relation on $X$ into which $\Ezero$ does not continuously embed, $F$ is
  a \Borel equivalence relation on $X$, and every $E$-class is a countable
  union of $(E \intersection F)$-classes. Then $X$ is a countable union of
  $(E \intersection F)$-invariant \Borel partial transversals of $E$ over $E
  \intersection F$.
\end{theorem}

\begin{theoremproof}
  Note that a set $Y \subseteq X$ is a partial transversal of $E$ over $E
  \intersection F$ if and only if it is independent with respect to the digraph
  $G = E \setminus F$. Moreover, Proposition \ref
  {generalization:separation} ensures that every $G$-independent \Borel
  set is contained in an $(E \intersection F)$-invariant $G$-independent
  \Borel set, so we need only show that there is a \Borel $\N$-coloring of
  $G$.
  
  Suppose, towards a contradiction, that there is no such coloring. Then
  Theorem \ref{gzero:main} yields a continuous homomorphism $\phi \from
  \Cantorspace \to X$ from $\Gzero$ to $G$. As the sets $E' = \preimage
  {(\phi \times \phi)}{E}$ and $F' = \preimage{(\phi \times \phi)}{F}$ have
  the \Baire property (see, for example, \cite[Theorem 21.6]{Kechris}),
  Proposition \ref{generalization:meager} ensures that $E'$ is meager. As
  the closed relation $D' = \preimage{(\phi \times \phi)}{\diagonal{X}}$ is
  contained in $E'$, it is nowhere dense, so Proposition \ref
  {generalization:homomorphism} yields a continuous homomorphism $\psi
  \from \Cantorspace \to \Cantorspace$ from $\triple{\setcomplement
  {\diagonal{\Cantorspace}}}{\Ezero}{\setcomplement{\Ezero}}$ to $\triple
  {\setcomplement{D'}}{\Ezero}{\setcomplement{E'}}$. As $\Gzero
  \subseteq E'$, Proposition \ref{generalization:generates} ensures that
  $\Ezero \subseteq E'$, so $\phi \composition \psi$ is a continuous
  embedding of $\Ezero$ into $E$, the desired contradiction.
\end{theoremproof}

\begin{remark}
  Under the weaker assumption that $F$ is co-analytic, one can run the
  same argument without Proposition \ref{generalization:separation} to
  obtain the weaker conclusion in which the \Borel partial transversals of
  $E$ over $E \intersection F$ need not be $(E \intersection F)$-invariant.
\end{remark}

A partial function $\phi \from X \partialto Y$ \definedterm{uniformizes} a set
$R \subseteq X \times Y$ if $\domain{\phi} = \image{\projection[X]}{R}$ and
$\graph{\phi} \subseteq R$.

\begin{proposition} \label{generalization:uniformization:single}
  Suppose that $X$ and $Y$ are \Polish spaces, $F$ is a strongly idealistic
  \Borel equivalence relation on $Y$, and $P \subseteq X \times Y$ is a
  \Borel set whose non-empty vertical sections are $F$-classes. Then
  $\image{\projection[X]}{P}$ is \Borel and there is a \Borel function $\phi
  \from \image{\projection[X]}{P} \to Y$ that uniformizes $P$.
\end{proposition}

\begin{propositionproof}
  Fix a witness $y \mapsto \calJ_y$ to the strong idealisticity of $F$, and
  define $\calI_x = \calJ_y$ for all $\pair{x}{y} \in P$. Given a \Borel set $R
  \subseteq P$, observe that the set $R' = \set{\pair{\pair{y}{x}}{z} \in (Y
  \times X) \times Y}[x \mathrel{P} y \mathand x \mathrel{R} z]$ is \Borel, so
  the set \heightcorrection{$S = \set{x \in X}[\exists y \in \verticalsection{P}
  {x} \ \verticalsection{R'}{\pair{y}{x}} \in \calJ_y]$} is analytic, the set
  \heightcorrection{$T = \set{x \in X}[\forall y \in \verticalsection{P}{x}
  \ \verticalsection{R'}{\pair{y}{x}} \in \calJ_y]$} is co-analytic, and
  $\verticalsection{R}{x} \in \calI_x \iff x \in S \iff x \in T$ for all $x \in \image
  {\projection[X]}{P}$. As noted in \cite[Theorem 18.6*]{Kechris:Corrections},
  the proof of \cite[Theorem 18.6]{Kechris} can therefore be easily modified
  to ensure that $\image{\projection[X]}{P}$ is \Borel, as well as to obtain
  a \Borel function $\phi \from \image{\projection[X]}{P} \to Y$ that
  uniformizes $P$.
\end{propositionproof}

We say that a partial function $T \from X / E \partialto X / E$ is \definedterm
{\Borel} if preimages of \Borel subsets of $X / E$ are \Borel.

\begin{proposition} \label{generalization:Borel}
  Suppose that $X$ is a \Polish space, $E$ is a strongly idealistic \Borel
  equivalence relation on $X$, and $T \from X / E \partialto X / E$ is a
  strongly \Borel partial function. Then $T$ is \Borel.
\end{proposition}

\begin{propositionproof}
  By Proposition \ref{generalization:uniformization:single}, the set $D =
  \union[\domain{T}]$ is \Borel and there is a \Borel function
  \heightcorrection{$\tilde{T} \from D \to X$} for which \heightcorrection
  {$\graph{T} = \saturation{\graph{\tilde{T}}}{\diagonal{X} \times E}$}. It
  therefore only remains to note that if $B \subseteq X / E$ is \Borel, then
  \heightcorrection{$\union[\preimage{T}{B}] = \preimage{\tilde{T}}{\union
  [B]}$}, so $\preimage{T}{B}$ is \Borel, thus $T$ is \Borel.
\end{propositionproof}

The \definedterm{composition} of partial functions $S \from Y \partialto Z$
and $T \from X \partialto Y$ is given by $(S \composition T)(x) = z \iff \exists
y \in Y \ (T(x) = y \mathand S(y) = z)$.

\begin{proposition} \label{generalization:composition}
  Suppose that $X$ is a \Polish space, $E$ is a strongly idealistic \Borel
  equivalence relation on $X$, and $S, T \from X / E \partialto X / E$ are
  strongly \Borel partial functions. Then $S \composition T$ is strongly
  \Borel.
\end{proposition}

\begin{propositionproof}
  By Proposition \ref{generalization:uniformization:single}, the sets $C =
  \union[\domain{S}]$ and $D = \union[\domain{T}]$ are \Borel and there
  are \Borel functions $\tilde{S} \from C \to X$ and $\tilde{T} \from D \to X$
  with $\graph{S} = \saturation{\graph{\tilde{S}}}{\diagonal{X} \times E}$ and
  $\graph{T} = \saturation{\graph{\tilde{T}}}{\diagonal{X} \times E}$. It only
  remains to note that if $x, y \in X$, then $(S \composition T)
  (\equivalenceclass{x}{E}) = \equivalenceclass{y}{E} \iff (x \in D
  \intersection \preimage{\tilde{T}}{C} \mathand (\tilde{S} \composition
  \tilde{T})(x) \mathrel{E} y)$, so $S \composition T$ is strongly \Borel.
\end{propositionproof}

The \definedterm{powers} of a partial injection $T \from X \partialto X$ are
given by $T^0 = \id[X]$, $T^{n+1} = T \composition T^n$ for all $n \in \N$,
and $T^{-n} = \inverse{(T^n)}$ for all $n \in \positiveintegers$. The
\definedterm{$T$-saturation} of a set $Y \subseteq X$ is given by
\heightcorrection{$\saturation{Y}{T} = \union[n \in \Z][\image{T^n}{\domain
{T^n} \intersection Y}]$}. The \definedterm{$T$-orbit} of a point $x \in X$ is
given by $\orbit{x}{T} = \saturation{\set{x}}{T}$ and the \definedterm{orbit
equivalence relation} induced by $T$ is the equivalence relation on $X$
given by $x \mathrel{\orbitequivalencerelation{T}{X}} y \iff \orbit{x}{T} =
\orbit{y}{T}$.

\begin{proposition} \label{generalization:bijection}
  Suppose that $X$ is a \Polish space, $E$ is a strongly idealistic \Borel
  equivalence relation on $X$, and $S \from X / E \partialto X / E$ is a
  strongly \Borel partial injection. Then there is a strongly \Borel bijection
  $T \from X / E \to X / E$ for which $\orbitequivalencerelation{S}{X / E} =
  \orbitequivalencerelation{T}{X / E}$.
\end{proposition}

\begin{propositionproof}
  Set $Y = \setcomplement{\domain{\inverse{S}}}$ and $Z =
  \setcomplement{\domain{S}}$, fix transitive permutations $\sigma$ of
  $\N$ and $\tau$ of $-\N$, and define $T \from X / E \to X / E$ by
  \begin{equation*}
    T(x) =
      \begin{cases}
        S(x) & \text{if $x \in \setcomplement{(\saturation{Y}{S} \union
          \saturation{Z}{S})} \union (\saturation{Y}{S} \intersection \saturation
            {Z}{S} \setminus Z)$,} \\
        \image{S^{-n}}{x} & \text{if $n \in \N$ and $x \in \image{S^n}{Y}
          \intersection Z$,} \\
        \image{S^{\sigma(n) - n}}{x} & \text{if $n \in \N$ and $x \in \image{S^n}
          {Y} \setminus \saturation{Z}{S}$, and} \\
        \image{S^{\tau(n) - n}}{x} & \text{if $n \in -\N$ and $x \in \image
          {S^n}{Z} \setminus \saturation{Y}{S}$.}
      \end{cases}
  \end{equation*}
  Propositions \ref{generalization:Borel} and \ref{generalization:composition}
  ensure that $T$ is strongly \Borel. The first clause in the definition of $T$
  ensures that the orbit equivalence relations of $S$ and $T$ agree on the
  set of points whose $S$-orbit has type $\Z$, the first two clauses ensure
  that the orbit equivalence relations agree on the set of points whose
  $S$-orbit is finite, the third clause ensures that the orbit equivalence
  relations agree on the set of points whose $S$-orbit has type $\N$, and
  the fourth clause ensures that the orbit equivalence relations agree on the
  set of points whose $S$-orbit has type $-\N$.
\end{propositionproof}

A \definedterm{reduction} of a binary relation $R$ on a set $X$ to a
binary relation $S$ on a set $Y$ is a homomorphism from $\pair{R}
{\setcomplement{R}}$ to $\pair{S}{\setcomplement{S}}$.

\begin{theorem} \label{generalization:uniformization}
  Suppose that $X$ and $Y$ are \Polish spaces, $F$ is a strongly idealistic
  \Borel equivalence relation on $Y$, and $R \subseteq X \times Y$ is a
  \Borel set whose vertical sections are countable unions of $F$-classes.
  Then $\image{\projection[X]}{R}$ is \Borel and there are \Borel maps
  $\phi_n \from \image{\projection[X]}{R} \to Y$ with the property that
  $\forall x \in \image{\projection[X]}{R} \ \verticalsection{R}{x} = \union[n \in
  \N][\equivalenceclass{\phi_n(x)}{F}]$.
\end{theorem}

\begin{theoremproof}
  Note that there is no continuous embedding $\pi \from \Cantorspace \into
  X \times Y$ of $\Ezero$ into $\diagonal{X} \times \completeer{Y}$, since
  otherwise $\projection[X] \composition \pi$ would be a continuous
  reduction of $\Ezero$ to $\diagonal{X}$, contradicting the well-known fact
  that every continuous homomorphism from $\Ezero$ to $\diagonal{X}$ is
  constant. Theorem \ref{generalization:transversal} therefore yields
  $(\diagonal{X} \times F)$-invariant \Borel partial transversals $P_n
  \subseteq R$ of $\diagonal{X} \times \completeer{Y}$ over $\diagonal{X}
  \times F$ for which $R = \union[n \in \N][P_n]$. For all $n \in \N$,
  Proposition \ref{generalization:uniformization:single} ensures that the set
  $D_n = \image{\projection[X]}{P_n}$ is \Borel and yields a \Borel function
  $\psi_n \from D_n \to Y$ that uniformizes $P_n$. Then $\image{\projection
  [X]}{R} = \union[n \in \N][D_n]$, so $\image{\projection[X]}{R}$ is \Borel,
  and the functions of the form $\phi_n = \psi_n \union \union[k \in \N]
  [{\restriction{\psi_k}{(D_k \setminus \union[j \in k \union \set{n}][D_j])}}]$
  are as desired.
\end{theoremproof}

\begin{remark}
  For all $n \in \N$, the set $R_n = \saturation{\graph{\phi_n}}{\diagonal{X}
  \times F} = \set{\pair{x}{y} \in \image{\projection[X]}{R} \times Y}[\phi_n(x)
  \mathrel{F} y]$ is \Borel, so Theorem \ref{generalization:uniformization}
  yields the special case of Theorem \ref{intro:uniformization} referred to in
  Remark \ref{intro:uniformization:E}.
\end{remark}

\begin{remark}
  Under the weaker assumption that $X$ and $Y$ are \Hausdorff, $F$ is
  co-analytic, and $R$ is analytic, the above argument still yields a
  conclusion similar to that of Theorem \ref{intro:uniformization}, although
  the resulting sets neither have full projections nor enjoy the same level of
  definability.
\end{remark}

We say that a subequivalence relation $F$ of an equivalence relation $E$
has \definedterm{countable index} if every $E$-class is a countable union
of $F$-classes.

\begin{proposition} \label{generalization:saturation}
  Suppose that $X$ is a \Polish space, $E$ is a \Borel equivalence relation
  on $X$, $F$ is a countable-index strongly idealistic \Borel subequivalence
  relation of $E$, and $B \subseteq X$ is an $F$-invariant \Borel set. Then
  $\saturation{B}{E}$ is \Borel.
\end{proposition}

\begin{propositionproof}
  By Theorem \ref{generalization:uniformization}, there are \Borel
  functions $\phi_n \from X \to X$ for which $E = \union[n \in \N][\saturation
  {\graph{\phi_n}}{\diagonal{X} \times F}]$. Then $\saturation{B}{E} = \union
  [n \in \N][\preimage{\phi_n}{B}]$.
\end{propositionproof}

Given an equivalence relation $E$ on a set $X$, we say that a set $Y
\subseteq X$ is \definedterm{$E$-complete} if $X = \saturation{Y}{E}$.
A \definedterm{transversal} of an equivalence relation $E$ over a
subequivalence relation $F$ is an $E$-complete partial transversal of $E$
over $F$.

\begin{proposition} \label{generalization:transversal:two}
  Suppose that $X$ is a topological space, $F \subseteq E$ are
  equivalence relations on $X$ for which the $E$-saturation of every
  $F$-invariant \Borel set is \Borel, there is a cover $\sequence{B_n}[n \in
  \N]$ of $X$ by $F$-invariant \Borel partial transversals of $E$ over $F$,
  and $B \subseteq X$ is an $F$-invariant \Borel partial transversal of $E$
  over $F$. Then $B$ is contained in an $F$-invariant \Borel transversal of
  $E$ over $F$.
\end{proposition}

\begin{propositionproof}
  Set $B_0' = B$, recursively define $B_{n+1}' = B_n' \union (B_n \setminus
  \saturation{B_n'}{E})$ for all $n \in \N$, and observe that the set $B' =
  \union[n \in \N][B_n']$ is as desired.
\end{propositionproof}

\section{Finite bases} \label{finitebasis}

The \definedterm{support} of a sequence $s \in \Cantortree$ is given
by $\support{s} = \preimage{s}{\set{1}}$.

\begin{proposition} \label{finitebasis:homomorphism}
  Suppose that $X$ is a set, $T \from X \to X$ is a bijection, and $\phi
  \from \Cantorspace \to X$ is a homomorphism from $\Gzero$ to
  the graph of $T$. Then $\forall n \in \N \forall u, v \in \Cantorspace[n]
  \forall c \in \Cantorspace \ T^{\cardinality{\support{v}}-\cardinality
  {\support{u}}}(\phi(u \concatenation c)) = \phi(v \concatenation c)$.
\end{proposition}

\begin{propositionproof}
  Suppose that we have already established the proposition at some $n \in
  \N$. To see that it holds at $n + 1$, observe that if $c \in \Cantorspace$
  and $u, v \in \Cantorspace[n]$, then $T^{\cardinality{\support{\bbs[n]}} -
  \cardinality{\support{u}}}(\phi(u \concatenation \sequence{0}
  \concatenation c)) = \phi(\bbs[n] \concatenation \sequence{0}
  \concatenation c)$, $T(\phi(\bbs[n] \concatenation \sequence{0}
  \concatenation c)) = \phi(\bbs[n] \concatenation \sequence{1}
  \concatenation c)$ since $\phi$ is a homomorphism, and $T^{\cardinality
  {\support{v}} - \cardinality{\support{\bbs[n]}}}(\phi(\bbs[n] \concatenation
  \sequence{1} \concatenation c)) = \phi(v \concatenation \sequence{1}
  \concatenation c)$. But $\cardinality{\support{\bbs[n]}} - \cardinality
  {\support{u}} + 1 + \cardinality{\support{v}} - \cardinality{\support{\bbs[n]}}
  = \cardinality{\support{v \concatenation \sequence{1}}} - \cardinality
  {\support{u \concatenation \sequence{0}}}$, so $T^{\cardinality{\support{v
  \concatenation \sequence{1}}} - \cardinality{\support{u \concatenation
  \sequence{0}}}}(\phi(u \concatenation \sequence{0} \concatenation c)) =
  \phi(v \concatenation \sequence{1} \concatenation c)$.
\end{propositionproof}

An equivalence relation $E$ on a space $X$ is \definedterm
{generically ergodic} if every $E$-invariant set with the \Baire property
is meager or comeager.

\begin{proposition} \label{finitebasis:genericallyergodic}
  Suppose that $k \ge 2$. Then $\Fzero{k}$ is generically ergodic.
\end{proposition}

\begin{propositionproof}
  Suppose that $B \subseteq \Cantorspace$ is a non-meager $\Fzero
  {k}$-invariant set with the \Baire property, fix a sequence $s \in
  \Cantortree$ for which $B$ is comeager in $\extensions{s}$, and set
  $n = \length{s}$. It is sufficient to show that $B$ is comeager in
  $\extensions{u}$ for all $u \in \Cantorspace[k-1+n]$. Towards this
  end, fix an extension $t \in \Cantorspace[k-1+n]$ of $s$ for which
  $\sum_{i < k-1+n} t(i) \equiv \sum_{i < k-1+n} u(i) \mod{k}$ and define
  $\iota \from \extensions{t} \to \extensions{u}$ by $\iota(t \concatenation c)
  = u \concatenation c$ for all $c \in \Cantorspace$. As $\iota$ is a
  homeomorphism, it follows that $\image{\iota}{B}$ is comeager in $\image
  {\iota}{\extensions{t}}$. But the former set is contained in $B$, and the
  latter is $\extensions{u}$.
\end{propositionproof}

\begin{proposition} \label{finitebasis:descriptivelyergodic}
  Suppose that $X$ is a \Baire space, $E$ is an equivalence relation
  on $X$ with respect to which saturations of meager sets are meager,
  $F$ is a generically-ergodic subequivalence relation of $E$, and $B
  \subseteq X$ is an $F$-invariant set with the \Baire property whose
  complement is $E$-complete. Then $B$ is meager.
\end{proposition}

\begin{propositionproof}
  As $\saturation{\setcomplement{B}}{E} = X$, it follows that $\saturation
  {\setcomplement{B}}{E}$ is not meager, so $\setcomplement{B}$ is not
  meager, thus $\setcomplement{B}$ is comeager.
\end{propositionproof}

We next establish a variant of Proposition \ref{generalization:saturation}:

\begin{proposition} \label{finitebasis:saturations}
  Suppose that $n \in \positiveintegers$, $X$ is a \Hausdorff space, $E$ is
  an analytic equivalence relation on $X$, $F$ is a co-analytic equivalence
  relation on $X$, $E$ has index $n$ over $E \intersection F$, and $B
  \subseteq X$ is an $(E \intersection F)$-invariant \Borel set. Then
  $\saturation{B}{E}$ is \Borel.
\end{proposition}

\begin{propositionproof}
  As saturations of analytic sets with respect to analytic equivalence
  relations are clearly analytic, it is sufficient to show that if $C \subseteq X$
  is an $(E \intersection F)$-invariant co-analytic set, then $\saturation{C}
  {E}$ is co-analytic. Towards this end, let $R$ be the set of pairs $\pair{x}
  {y} \in X \times \functions{n}{X}$ such that $\forall i < n \ x \mathrel{E}
  y(i)$, $\forall i < j < n \ \neg y(i) \mathrel{F} y(j)$, and $\forall i < n \ y(i)
  \notin C$, and note that $\setcomplement{\saturation{C}{E}} = \image
  {\projection[X]}{R}$ and $R$ is analytic, thus so too is $\setcomplement
  {\saturation{C}{E}}$.
\end{propositionproof}

A \definedterm{reduction} of a sequence $\sequence{R_i}[i \in I]$
of binary relations on a set $X$ to a sequence $\sequence{S_i}[i \in I]$ of
binary relations on a set $Y$ is a function $\phi \from X \to Y$ that is a
reduction of $R_i$ to $S_i$ for all $i \in I$. An \definedterm{embedding}
is an injective reduction. We next note that the family $\set{\pair{\Ezero}
{\Fzero{p}}}[p \text{ is prime}]$ is a minimal basis for $\set{\pair{\Ezero}
{\Fzero{k}}}[k \ge 2]$ under any quasi-order between \Baire-measurable
reducibility and continuous embeddability:

\begin{proposition} \label{finitebasis:examples}
  Suppose that $j, k \ge 2$.
  \begin{enumerate}
    \item If $j \divides k$, then there is a continuous embedding
      $\phi \from \Cantorspace \to \Cantorspace$ of $\pair{\Ezero}
      {\Fzero{j}}$ into $\pair{\Ezero}{\Fzero{k}}$.
    \item If $j$ and $k$ are relatively prime, then there is no
      \Baire-measurable homomorphism $\phi \from \Cantorspace \to
      \Cantorspace$ from $\pair{\Ezero \setminus \Fzero{j}}{\Fzero{j}}$
      to $\pair{\Ezero \setminus \Fzero{k}}{\Fzero{k}}$.
  \end{enumerate}
\end{proposition}

\begin{propositionproof}
  To see (1), observe that the function $\phi \from \Cantorspace \to
  \Cantorspace$ given by $\phi(c) = \concatenation[n \in \N]
  [\constantsequence{c(n)}{k / j}]$ is a continuous embedding of $\pair
  {\Ezero}{\Fzero{j}}$ into $\pair{\Ezero}{\Fzero{k}}$.
  
  To see (2), suppose, towards a contradiction, that there is such a
  homomorphism $\phi \from \Cantorspace \to \Cantorspace$. For
  all $\ell \in \set{j, k}$, define $\psi_\ell \from \Cantorspace \to
  \Cantorspace$ by $\psi_\ell(\constantsequence{1}{n} \concatenation
  \sequence{0} \concatenation c) = \constantsequence{1}{n}
  \concatenation \sequence{1} \concatenation c$ for all $c \in
  \Cantorspace$ and $n \in \N$ and $\psi_\ell(\constantsequence{1}
  {\infty}) = \constantsequence{0}{\ell-1} \concatenation
  \constantsequence{1}{\infty}$, so that $\psi_\ell$ factors over $\Fzero
  {\ell}$ to a permutation of $\Cantorspace / \Fzero{\ell}$ whose orbit
  equivalence relation is $\Ezero / \Fzero{\ell}$. In particular, it follows
  that the sets of the form $B_i = \set{c \in \Cantorspace}[\phi(\psi_j(c))
  \mathrel{\Fzero{k}} \psi_k^i(\phi(c))]$, for $0 < i < k$, are $\Fzero
  {j}$-invariant and partition $\Cantorspace$, so Proposition \ref
  {finitebasis:genericallyergodic} yields a unique $i$ for which $B_i$ is
  comeager. As $\Ezero$-saturations of meager sets are meager, there
  are comeagerly-many $c \in \Cantorspace$ for which $\equivalenceclass
  {c}{\Ezero} \subseteq B_i$. Given any such $c$, note that $c \mathrel
  {\Fzero{j}} \psi_j^j(c)$, so $\phi(c) \mathrel{\Fzero{k}} \psi_k^{ij}(\phi(c))$,
  thus $k \divides ij$, the desired contradiction.
\end{propositionproof}

Along with Propositions \ref{generalization:saturation} and \ref
{generalization:transversal:two}, the following fact yields the special case
of Theorem \ref{intro:transversal} where there exists $k \in \N$ for which
every $E$-class has cardinality at most $k$:

\begin{theorem} \label{finitebasis:transversal}
  Suppose that $k \ge 2$, $X$ is a \Hausdorff space, $E$ is an analytic
  equivalence relation on $X$, $F$ is a \Borel equivalence relation on $X$,
  and every $E$-class is a union of at most $k$ $(E \intersection
  F)$-classes. Then exactly one of the following holds:
  \begin{enumerate}
    \item There is a cover $\sequence{B_i}[i < k]$ of $X$ by $(E \intersection
      F)$-invariant \Borel partial transversals of $E$ over $E \intersection F$.
    \item There is a continuous embedding $\pi \from \Cantorspace \into X$
      of $\pair{\Ezero}{\Fzero{p}}$ into $\pair{E}{F}$ for some prime $p \le k$.
  \end{enumerate}
\end{theorem}

\begin{theoremproof}
  To see that the conditions are mutually exclusive, note that if both hold,
  then there exists $i < k$ for which $\preimage{\pi}{B_i}$ is a non-meager
  \Borel partial transversal of $\Ezero$ over $\Fzero{p}$, contradicting
  Proposition \ref{finitebasis:descriptivelyergodic}.

  To see that at least one of the conditions holds, suppose that we have
  already established the theorem strictly below $k$, and define $Y = \set{y
  \in \functions{k}{X}}[\forall i < j < k \ y(i) \mathrel{(E \setminus F)} y(j)]$. It
  is sufficient to show that at least one of the following holds:

  \begin{enumerate}
    \renewcommand{\theenumi}{\alph{enumi}}
    \item There is an $(E \intersection F)$-invariant \Borel partial transversal
      $B \subseteq X$ of $E$ over $E \intersection F$ for which $\image
      {\projection[0]}{Y} \subseteq \saturation{B}{E}$.
    \item There is a continuous embedding $\pi \from \Cantorspace \into
      \image{\projection[0]}{Y}$ of $\pair{\Ezero}{\Fzero{p}}$ into $\pair
      {\restriction{E}{\image{\projection[0]}{Y}}}{\restriction{F}{\image
      {\projection[0]}{Y}}}$ for some prime $p \le k$.
  \end{enumerate}
  
  Let $\Sigma$ be the set of permutations $\sigma$ of $k$ for which
  $\sigma(0) \neq 0$, and for each $\sigma \in \Sigma$, let $G_\sigma$
  be the digraph on $Y$ with respect to which two sequences $y$ and
  $z$ are related if and only if $\forall i < k \ y(\sigma(i)) \mathrel{(E
  \intersection F)} z(i)$.
  
  \begin{lemma}
    Suppose that there are \Borel $\N$-colorings $c_\sigma \from Y \to
    \N$ of $G_\sigma$ for all $\sigma \in \Sigma$. Then there is an $(E
    \intersection F)$-invariant \Borel partial transversal $B \subseteq X$
    of $E$ over $E \intersection F$ for which $\image{\projection[0]}{Y}
    \subseteq \saturation{B}{E}$.
  \end{lemma}
  
  \begin{lemmaproof}
    Note that the function $c_\Sigma \from Y \to \functions{\Sigma}{\N}$
    given by $c_{\Sigma}(y)(\sigma) = c_\sigma(y)$ is a \Borel $\functions
    {\Sigma}{\N}$-coloring of the digraph $G_\Sigma = \union[\sigma \in
    \Sigma][G_\sigma]$, and if $Z \subseteq Y$ is a \Borel
    $G_\Sigma$-independent set, then $\image{\projection[0]}{Z}$ is an
    analytic partial transversal of $E$ over $E \intersection F$, so
    Proposition \ref{generalization:separation} ensures that it is contained in
    an $(E \intersection F)$-invariant \Borel partial transversal of $E$ over
    $E \intersection F$, thus there is a cover $\sequence{B_n}[n \in \N]$ of
    $\image{\projection[0]}{Y}$ by $(E \intersection F)$-invariant \Borel
    partial transversals of $E$ over $E \intersection F$. Then Propositions
    \ref{generalization:transversal:two} and \ref{finitebasis:saturations} yield
    an $(E \intersection F)$-invariant \Borel transversal $C \subseteq \image
    {\projection[0]}{Y}$ of $\restriction{E}{\image{\projection[0]}{Y}}$ over
    $\restriction{(E \intersection F)}{\image{\projection[0]}{Y}}$, in which
    case one more application of Proposition \ref{generalization:separation}
    yields an $(E \intersection F)$-invariant \Borel partial transversal $B
    \supseteq C$ of $E$ over $E \intersection F$.
  \end{lemmaproof}
  
  Suppose now that there exists $\sigma \in \Sigma$ for which there is no
  \Borel coloring $c_\sigma \from Y \to \N$ of $G_\sigma$. Then Theorem
  \ref{gzero:main} yields a continuous homomorphism $\phi \from
  \Cantorspace \to Y$ from $\Gzero$ to $G_\sigma$. Define $\phi_0 =
  \projection[0] \composition \phi$ and $j = \cardinality{\set{\sigma^n(0)}[n
  \in \Z]}$. As the relations $E' = \preimage{(\phi_0 \times \phi_0)}{E}$ and
  $F' = \preimage{(\phi_0 \times \phi_0)}{F}$ have the \Baire property,
  Proposition \ref{generalization:meager} ensures that they are meager, and
  since the closed relation $D' = \preimage{(\phi_0 \times \phi_0)}{\diagonal
  {X}}$ is contained in $E'$, it is nowhere dense, so Proposition \ref
  {generalization:homomorphism} yields a continuous homomorphism $\psi
  \from \Cantorspace \to \Cantorspace$ from $\quadruple{\setcomplement
  {\diagonal{\Cantorspace}}}{\Fzero{j}}{\Ezero \setminus \Fzero{j}}
  {\setcomplement{\Ezero}}$ to $\quadruple{\setcomplement{D'}}{\Fzero{j}}
  {\Ezero \setminus \Fzero{j}}{\setcomplement{(E' \union F')}}$. As
  Proposition \ref{finitebasis:homomorphism} ensures that $\phi_0$ is a
  homomorphism from $\pair{\Fzero{j}}{\Ezero \setminus \Fzero{j}}$ to
  $\pair{E \intersection F}{E \setminus F}$, it follows that $\phi_0
  \composition \psi$ is a continuous embedding of $\pair{\Ezero}{\Fzero{j}}$
  into $\pair{E}{F}$, so it only remains to observe that if $p$ is any prime
  dividing $j$, then Proposition \ref{finitebasis:examples} ensures that there
  is a continuous embedding of $\pair{\Ezero}{\Fzero{p}}$ into $\pair
  {\Ezero}{\Fzero{j}}$, and therefore of $\pair{\Ezero}{\Fzero{p}}$ into $\pair
  {E}{F}$.
\end{theoremproof}

A \definedterm{partial uniformization} of a set $R \subseteq X \times Y$
over an equivalence relation $F$ on $Y$ is a subset of $R$ whose vertical
sections are contained in $F$-classes. Along with Theorem \ref
{generalization:uniformization} and Proposition \ref
{generalization:transversal:two}, the following fact yields the special case
of Theorem \ref{intro:uniformization} where there exists $k \in \N$ for which
every vertical section of $R$ has cardinality at most $k$:

\begin{theorem} \label{finitebasis:uniformization}
  Suppose that $k \ge 2$, $X$ and $Y$ are \Hausdorff spaces, $E$ is an
  analytic equivalence relation on $X$, $F$ is a \Borel equivalence relation
  on $Y$, and $R \subseteq X \times Y$ is an $(E \times \diagonal
  {Y})$-invariant analytic set whose vertical sections are contained in unions
  of at most $k$ $F$-classes. Then exactly one of the following holds:
  \begin{enumerate}
    \item There is a cover $\sequence{R_i}[i < k]$ of $R$ by $(\restriction{(E
      \times F)}{R})$-invariant \Borel-in-$R$ partial uniformizations of $R$
      over $F$.
    \item There are continuous embeddings $\pi_X \from \Cantorspace \into
      X$ of $\Ezero$ into $E$ and $\pi_Y \from \Cantorspace \into Y$ of
      $\Fzero{p}$ into $F$ such that $\image{(\pi_X \times \pi_Y)}{\Ezero}
      \subseteq R$ for some prime $p \le k$.
  \end{enumerate}
\end{theorem}

\begin{theoremproof}
  To see that the conditions are mutually exclusive, note that if both hold,
  then there exists $i < k$ for which $\image{\projection[1]}{\preimage
  {(\pi_X \times \pi_Y)}{R_i}}$ is a non-meager partial transversal of
  $\Ezero$ over $\Fzero{p}$ with the \Baire property, contradicting
  Proposition \ref{finitebasis:descriptivelyergodic}.
  
  To see that at least one of the conditions holds, note that
  $\restriction{(E \times \completeer{Y})}{R}$ is analytic and $\restriction
  {(\completeer{X} \times F)}{R}$ is \Borel, and appeal to Theorem \ref
  {finitebasis:transversal} to see that if condition (1) fails, then there is a
  continuous embedding $\pi \from \Cantorspace \into R$ of $\pair{\Ezero}
  {\Fzero{p}}$ into $\pair{\restriction{(E \times \completeer{Y})}{R}}
  {\restriction{(E \times F)}{R}}$ for some prime $p \le k$. It follows that the
  function $\pi_X' = \projection[X] \composition \pi$ is a continuous
  reduction of $\Ezero$ to $E$ and the function $\pi_Y' = \projection[Y]
  \composition \pi$ is a continuous homomorphism from $\pair{\Fzero{p}}
  {\Ezero \setminus \Fzero{p}}$ to $\pair{F}{\setcomplement{F}}$, and
  therefore from $\Gzero$ to $\setcomplement{F}$, so Proposition \ref
  {generalization:meager} ensures that the equivalence relation $F' =
  \preimage{(\pi_Y' \times \pi_Y')}{F}$ is meager, in which case the closed
  relations $D_X' = \preimage{(\pi_X' \times \pi_X')}{\diagonal{X}}$ and
  $D_Y' = \preimage{(\pi_Y' \times \pi_Y')}{\diagonal{Y}}$ are nowhere
  dense, so Proposition \ref{generalization:homomorphism} yields a
  continuous homomorphism $\pi' \from \Cantorspace \to \Cantorspace$
  from $\quadruple{\setcomplement{\diagonal{\Cantorspace}}}{\Fzero{p}}
  {\Ezero \setminus \Fzero{p}}{\setcomplement{\Ezero}}$ to $\quadruple
  {\setcomplement{(D_X' \union D_Y')}}{\Fzero{p}}{\Ezero \setminus \Fzero
  {p}}{\setcomplement{(\Ezero \union F')}}$, thus the functions $\pi_X =
  \pi_X' \composition \pi'$ and $\pi_Y = \pi_Y' \composition \pi'$ are
  continuous embeddings of $\Ezero$ into $E$ and $\Fzero{p}$ into $F$.
  As $\image{\pi}{\Cantorspace} \subseteq R$, it follows that $\image{(\pi_X'
  \times \pi_Y')}{\diagonal{\Cantorspace}} \subseteq R$, so the facts that
  $\pi_X'$ is a homomorphism from $\Ezero$ to $E$ and $R$ is $(E \times
  \diagonal{Y})$-invariant ensure that $\image{(\pi_X' \times \pi_Y')}
  {\Ezero} \subseteq R$, thus the fact that $\pi'$ is a homomorphism from
  $\Ezero$ to $\Ezero$ implies that  $\image{(\pi_X \times \pi_Y)}{\Ezero}
  \subseteq R$.
\end{theoremproof}

Given a binary relation $R$ on $X$, set $\inverse{R} = \set{\pair{y}{x} \in X
\times X}[x \mathrel{R} y]$. Given an equivalence relation $F$ on $X$, we
say that $R$ is the \definedterm{graph of a partial injection} over $F$ if
$\forall \pair{x}{y}, \pair{x'}{y'} \in R \ ( x \mathrel{F} x' \iff y \mathrel{F} y')$.
Along with Proposition \ref{generalization:bijection}, the following fact
yields the special case of Theorem \ref{intro:graph} where there exists $k
\in \N$ for which every $E$-class has cardinality at most $k + 1$:

\begin{theorem} \label{finitebasis:graph}
  Suppose that $k \ge 2$, $X$ is a \Hausdorff space, $E$ is an analytic
  equivalence relation on $X$, $F$ is a \Borel equivalence relation on $X$,
  and every $E$-class is a union of at most $k + 1$ $(E \intersection
  F)$-classes. Then exactly one of the following holds:
  \begin{enumerate}
    \item There is a cover $\sequence{R_{i,j}}[i,j<k]$ of $E \setminus F$ by
      $((E \intersection F) \times (E \intersection F))$-invariant \Borel graphs
      of partial injections over $E \intersection F$.
    \item There is a continuous embedding $\pi \from \Cantorspace \times 2
      \into X$ of $\pair{\Ezero \times \completeer{2}}{\Ezero \disjointunion
      \Fzero{p}}$ into $\pair{E}{F}$ for some prime $p \le k$.
  \end{enumerate}
\end{theorem}

\begin{theoremproof}
  To see that the conditions are mutually exclusive, observe that if they both
  hold, then there exist $i, j < k$ with the property that $\set{c \in
  \Cantorspace}[\pi(c,0) \mathrel{R_{i,j}} \pi(c,1)] \times \set{1}$ is a
  non-meager partial transversal of $\Ezero \times \diagonal{\set{1}}$ over
  $\Fzero{p} \times \diagonal{\set{1}}$ with the \Baire property, a
  contradiction.

  To see that at least one of the conditions holds, note first that if there are
  $\restriction{((E \intersection F) \times F)}{(E \setminus F)}$-invariant
  \Borel-in-$(E \setminus F)$ partial uniformizations $R_j$ of $E \setminus
  F$ over $F$ for which \heightcorrection{$E \setminus F = \union[j < k]
  [R_j]$}, then the sets of the form \heightcorrection{$R_i \intersection
  \inverse{R_j}$}, where $i, j < k$, are as desired. By Theorem \ref
  {finitebasis:uniformization}, we can therefore assume that there are
  continuous embeddings $\pi_X \from \Cantorspace \into X$ of $\Ezero$
  into $E \intersection F$ and $\pi_Y \from \Cantorspace \into X$ of $\Fzero
  {p}$ into $F$ such that $\image{(\pi_X \times \pi_Y)}{\Ezero} \subseteq E
  \setminus F$ for some prime $p \le k$. As $\Ezero$ has countable index
  below the equivalence relation $E' = \preimage{(\pi_X \times \pi_X)}{E}$,
  the latter is meager, so the closed subequivalence relation $D' =
  \preimage{(\pi_X \times \pi_X)}{\diagonal{X}}$ is nowhere dense, thus
  Proposition \ref{generalization:homomorphism} yields a continuous
  homomorphism $\pi' \from \Cantorspace \to \Cantorspace$ from
  $\quadruple{\setcomplement{\diagonal{\Cantorspace}}}{\Fzero{p}}{\Ezero
  \setminus \Fzero{p}}{\setcomplement{\Ezero}}$ to $\quadruple
  {\setcomplement{D'}}{\Fzero{p}}{\Ezero \setminus \Fzero{p}}
  {\setcomplement{E'}}$. Define $\pi \from \Cantorspace \times 2 \to X$ by
  $\pi(c, 0) = (\pi_X \composition \pi')(c)$ and $\pi(c, 1) = (\pi_Y \composition
  \pi')(c)$ for all $c \in \Cantorspace$.
\end{theoremproof}

\section{Dichotomies} \label{dichotomy}

Let \definedsymbol{\sets{n}{X}} denote the family of subsets of $X$ of
cardinality $n$. When $X$ is a topological space, we endow $\sets{n}
{X}$ with the topology generated by the sets of the form $\set{a \in
\sets{n}{X}}[\exists f \from a \into n \forall x \in a \ x \in U_{f(x)}]$, where
$\sequence{U_i}[i < n]$ is a sequence of open subsets of $X$. Let
\definedsymbol{\sets{n}{X}[E]} denote the subspace consisting of
sets which are contained in a single $E$-class, and let \definedsymbol
{\sets{n}{X}[E][F]} denote the further subspace consisting of such sets
which are partial transversals of $F$. Define $\sets{<\aleph_0}{X} = \union
[n \in \N][\sets{n}{X}]$, $\sets{<\aleph_0}{X}[E] = \union[n \in \N][{\sets{n}{X}
[E]}]$, and $\sets{<\aleph_0}{X}[E][F] = \union[n \in \N][{\sets{n}{X}[E][F]}]$.

The \definedterm{trace} of a family $\calA \subseteq \sets{<\aleph_0}
{X}[E]$ on an equivalence class $C$ of $E$ is given by $\restriction
{\calA}{C} = \set{a \in \calA}[a \subseteq C]$. We say that two subsets
of $X$ are \definedterm{$F$-disjoint} if their $F$-saturations are
disjoint, $\calA$ is \definedterm{$F$-intersecting} if no two sets in
$\calA$ are $F$-disjoint, and $\calA$ is \definedterm{$E$-locally
$F$-intersecting} if its trace on every $E$-class is $F$-intersecting. A
set $Y \subseteq X$ \definedterm{punctures} $\calA$ if it intersects
every set in $\calA$. For each non-empty set $a' \in \sets{<\aleph_0}
{X}[E]$, define $\interval{a'}{\calA}[F] = \set{\saturation{a}{F}}[a \in
\calA \mathand a' \subseteq \saturation{a}{F}]$.

A \definedterm{partial quasi-transversal} of an equivalence relation $E$
over a subequivalence relation $F$ is a set $Y \subseteq X$ for which
there exists $k \in \N$ such that every $(\restriction{E}{Y})$-class is
contained in a union of at most $k$ $F$-classes. The following fact is
essentially a special case of \cite[Proposition 2.3.1]{CCM}; we provide the
proof for the reader's convenience:

\begin{proposition} \label{dichotomy:intersecting}
  Suppose that $X$ is a \Hausdorff space, $E$ is an analytic equivalence
  relation on $X$, $F$ is a \Borel equivalence relation on $X$, and $\calA
  \subseteq \sets{<\aleph_0}{X}[E][F]$ is an $E$-locally $F$-intersecting
  analytic family of sets of bounded finite cardinality. Then there is an $(E
  \intersection F)$-invariant Borel partial quasi-transversal $B \subseteq X$
  of $E$ over $E \intersection F$ puncturing $\calA$.
\end{proposition}

\begin{propositionproof}
  Recursively define $f(n) = n^{n+1} + \sum_{0 < k < n} f(k)$ for all $n \in
  \positiveintegers$. We will show that if every set in $\calA$ has cardinality
  at most $n$, then there is an $(E \intersection F)$-invariant \Borel set $B
  \subseteq X$ puncturing $\calA$ such that every $(\restriction{E}
  {B})$-class is contained in a union of at most $f(n)$ $(E \intersection
  F)$-classes.

  We proceed by induction on $n$. The base case $n = 1$ follows from an
  application of Proposition \ref{generalization:separation} to the
  equivalence relation $E \intersection F$ and the digraph $E \setminus F$,
  so suppose that $n \ge 2$ and we have established the proposition strictly
  below $n$. We will construct analytic families $\calA_k \subseteq \calA$
  and $\calA_k' \subseteq \sets{k}{X}[E][F]$, as well as $(E \intersection
  F)$-invariant \Borel sets $B_k \subseteq X$, such that:
  \begin{enumerate}
    \item $\forall k < n \ \calA_k = \set{a \in \calA_{k+1}}[a \intersection
      B_{k+1} = \emptyset]$.
    \item $\forall 1 \le k \le n \ \calA_k'  = \set{a' \in \sets{k}{X}[E][F]}
      [\cardinality{\interval{a'}{\calA_k}[E \intersection F]} > n^{n - k}]$.
    \item $\forall 1 \le k \le n \ B_k$ punctures $\calA_k'$.
    \item $\forall 1 \le k \le n \ E$ has index at most $f(k)$ over $E
      \intersection F$ on $B_k$.
  \end{enumerate}
  We proceed by reverse recursion, beginning with $\calA_n = \calA$,
  $\calA_n' = \emptyset$, and $B_n = \emptyset$. Suppose now that
  $0 < k < n$ and we have already found $\calA_{k+1}$, $\calA_{k+
  1}'$, and $B_{k+1}$. Conditions (1) and (2) then yield $\calA_k$ and
  $\calA_k'$.
  
  \begin{lemma} \label{dichotomy:intersecting:singleextension}
    Suppose that $a' \in \calA_k'$ and $x \in \saturation{a'}{E} \setminus
    \saturation{a'}{F}$. Then the family $\interval{a' \union \singleton{x}}
    {\calA_k}[E \intersection F]$ has cardinality at most $n^{n - (k + 1)}$.
  \end{lemma}
  
  \begin{lemmaproof}
    We can clearly assume that $\interval{a' \union \singleton{x}}
    {\calA_k}[E \intersection F] \neq \emptyset$. Given any set $a \in
    \interval{a' \union \singleton{x}}{\calA_k}[E \intersection F]$,
    condition (1) ensures that $a \intersection B_{k+1} =
    \emptyset$, thus $(a' \union \singleton{x}) \intersection B_{k+1} =
    \emptyset$. Condition (3) therefore implies that $a' \union
    \singleton{x} \notin \calA_{k+1}'$, so $\cardinality{\interval{a' \union
    \singleton{x}}{\calA_{k+1}}[E \intersection F]} \le n^{n - (k + 1)}$ by
    condition (2). As condition (1) also ensures that $\calA_k \subseteq
    \calA_{k+1}$, the lemma follows.
  \end{lemmaproof}
  
  \begin{lemma} \label{dichotomy:intersecting:setextension}
    Suppose that $a' \in \calA_k'$ and $b \subseteq \saturation{a'}{E}
    \setminus \saturation{a'}{F}$ has cardinality at most $n$. Then
    there exists $a \in \calA_k$ whose $(E \intersection F)$-saturation
    contains $a'$ and is disjoint from $b$.
  \end{lemma}
  
  \begin{lemmaproof}
    Lemma \ref{dichotomy:intersecting:singleextension} ensures that
    $\cardinality{\interval{a' \union \singleton{x}}{\calA_k}[E \intersection
    F]} \le n^{n - (k + 1)}$ for all $x \in b$, so $\interval{a'}{\calA_k}[E
    \intersection F]$ contains at most $n^{n - k}$ sets intersecting
    $b$, thus condition (2) implies that $\interval{a'}{\calA_k}[E
    \intersection F]$ contains a set disjoint from $b$.
  \end{lemmaproof}
  
  \begin{lemma}
    The family $\calA_k'$ is $E$-locally $F$-intersecting.
  \end{lemma}
  
  \begin{lemmaproof}
    Suppose, towards a contradiction, that there are $(E \intersection
    F)$-disjoint sets $a', b' \in \calA_k'$ with the property that
    $\saturation{a'}{E} = \saturation{b'}{E}$. Then Lemma \ref
    {dichotomy:intersecting:setextension} yields $b \in \calA_k$ whose
    $(E \intersection F)$-saturation contains $b'$ and is disjoint from
    $a'$, as well as $a \in \calA_k$ whose $(E \intersection
    F)$-saturation contains $a'$ and is disjoint from $b$, contradicting
    the fact that $\calA$ is $E$-locally $F$-intersecting.
  \end{lemmaproof}
  
  As the induction hypothesis yields an $(E \intersection F)$-invariant
  \Borel set $B_k \subseteq X$ puncturing $\calA_k'$ on which $E$
  has index at most $f(k)$ over $E \intersection F$, this completes the
  recursive construction. Define $A_0 = \union[\calA_0]$.
  
  \begin{lemma}
    Suppose that $x \in A_0$. Then $A_0 \intersection \equivalenceclass
    {x}{E}$ is contained in a union of at most $n^{n+1}$ $(E \intersection
    F)$-classes.
  \end{lemma}
  
  \begin{lemmaproof}
    Fix $a_0 \in \calA_0$ for which $x \in a_0$, and note that if $y \in
    a_0$, then $\cardinality{\set{\saturation{a}{E \intersection F}}[a \in
    \calA_0 \mathand y \in \saturation{a}{E \intersection F}]} \le n^{n-1}$,
    in which case $\union[{\set{\saturation{a}{E \intersection F}}[a \in
    \calA_0 \mathand y \in \saturation{a}{E \intersection F}]}]$ is a union
    of at most $n^n$ $(E \intersection F)$-classes. But every element of
    $A_0 \intersection \equivalenceclass{x}{E}$ is contained in a set of
    this form, so $A_0 \intersection \equivalenceclass{x}{E}$ is
    contained in a union of at most $n^{n+1}$ $(E \intersection
    F)$-classes.
      \end{lemmaproof}
  
  By applying Proposition \ref{generalization:separation} to the equivalence
  relation $E \intersection F$ and the dihypergraph $G = \set{x \in \functions
  {n^{n+1} + 1}{X}}[\forall i < j \le n^{n+1} \ x(i) \mathrel{(E \setminus F)}
  y(j)]$, we obtain an $(E \intersection F)$-invariant \Borel set $B_0
  \supseteq A_0$ with the property that every $(\restriction{E}{B_0})$-class
  is contained in a union of at most $n^{n+1}$ $(E \intersection F)$-classes.
  It only remains to note that if $a \in \calA$, then there is a least $k \le n$
  such that $a \in \calA_k$, in which case $a \intersection B_k$ is
  non-empty, so the set $B = \union[k \le n][B_k]$ punctures $\calA$.
\end{propositionproof}

A \definedterm{partial quasi-transversal} of an equivalence relation $E$ on
$X$ is a partial quasi-transversal of $E$ over $\diagonal{X}$. Theorem \ref
{intro:transversal} is a consequence of Propositions \ref
{generalization:saturation} and \ref{generalization:transversal:two},
Theorem \ref{finitebasis:transversal}, and the following:

\begin{theorem} \label{dichotomy:quasitransversal}
  Suppose that $X$ is a \Hausdorff space, $E$ is an analytic
  equivalence relation on $X$, $F$ is a \Borel equivalence relation on
  $X$, and every $E$-class is a countable union of $(E \intersection
  F)$-classes. Then exactly one of the following holds:
  \begin{enumerate}
    \item There is a cover $\sequence{B_n}[n \in \N]$ of $X$ by $(E
      \intersection F)$-invariant \Borel partial quasi-transversals of $E$ over
      $E \intersection F$.
    \item There exists a continuous embedding $\pi \from \Cantorspace
      \into X$ of $\pair{\Ezero}{\diagonal{\Cantorspace}}$ into $\pair
      {E}{F}$.
  \end{enumerate}
\end{theorem}

\begin{theoremproof}
  To see that the conditions are mutually exclusive, note that if both hold,
  then there exists $n \in \N$ for which $\preimage{\pi}{B_n}$ is a
  non-meager \Borel partial quasi-transversal of $\Ezero$. But the proof of
  the well-known fact that every partial transversal of $\Ezero$ with the
  \Baire property is meager works just as well to show that every partial
  quasi-transversal of $\Ezero$ with the \Baire property is meager, a
  contradiction.

  To see that at least one of the conditions holds, we can assume that $E
  \nsubseteq F$, so there are continuous surjections $\phi_{E \setminus F}
  \from \Bairespace \onto E \setminus F$ and $\phi_X \from \Bairespace
  \onto X$. We will recursively define a decreasing sequence $\sequence
  {B^\alpha}[\alpha < \omega_1]$ of \Borel subsets of $X$ whose
  complements are countable unions of $(E \intersection F)$-invariant
  \Borel partial quasi-transversals of $E$ over $E \intersection F$. We begin
  by setting $B^0 = X$. For all limit ordinals $\lambda < \omega_1$, we set
  $B^\lambda = \intersection[\alpha < \lambda][B^\alpha]$. To describe the
  construction at successor ordinals, we require several preliminaries.
  
  An \definedterm{approximation} is a triple of the form $a = \triple{n^a}
  {\phi^a}{\sequence{\psi_n^a}[n < n^a]}$, where $n^a \in \N$, $\phi^a
  \from \Cantorspace[n^a] \to \Bairespace[n^a]$, and $\psi_n^a \from
  \completeer{\Cantorspace[n]} \times \Cantorspace[n^a -
  (n+1)] \to \Bairespace[n^a]$ for all $n < n^a$. A \definedterm{one-step
  extension} of $a$ is an approximation $b$ for which:
  \begin{enumerate}
    \item $n^b = n^a + 1$.
    \item $\forall s \in \Cantorspace[n^a] \forall t \in \Cantorspace[n^b]
      \ (s \extendedby t \implies \phi^a(s) \extendedby
        \phi^b(t))$.
    \item $\forall n < n^a \forall r \in \completeer{\Cantorspace[n]}
      \forall s \in \Cantorspace[n^a - (n + 1)] \forall t \in \Cantorspace
        [n^b - (n + 1)]$ \\
          \hspace*{5pt} $(s \extendedby t \implies \psi^a_n(r, s)
            \extendedby \psi^b_n(r, t))$.
  \end{enumerate}
  A \definedterm{configuration} is a triple of the form $\gamma = \triple
  {n^\gamma}{\phi^\gamma}{\sequence{\psi_n^\gamma}[n <
  n^\gamma]}$, where $n^\gamma \in \N$, $\phi^\gamma \from
  \Cantorspace[n^\gamma] \to \Bairespace$, $\psi_n^\gamma \from
  \completeer{\Cantorspace[n]} \times \Cantorspace[n^\gamma - (n+1)]
  \to \Bairespace$ for all $n < n^\gamma$, and $(\phi_{E \setminus F}
  \composition \psi_n^\gamma)(r, s) = \sequence{(\phi_X \composition
  \phi^\gamma)(r(i) \concatenation \sequence{i} \concatenation s)}[i <
  2]$ for all $n < n^\gamma$, $r \in \completeer{\Cantorspace[n]}$, and
  $s \in \Cantorspace[n^\gamma - (n + 1)]$. We say that $\gamma$ is
  \definedterm{compatible} with a set $X' \subseteq X$ if $\image
  {(\phi_X \composition \phi^\gamma)}{\Cantorspace[n^\gamma]}
  \subseteq X'$, and \definedterm{compatible} with $a$ if:
  \begin{enumerate}
    \renewcommand{\theenumi}{\roman{enumi}}
    \item $n^a = n^\gamma$.
    \item $\forall s \in \Cantorspace[n^a] \ \phi^a(s) \extendedby
      \phi^\gamma(s)$.
    \item $\forall n < n^a \forall r \in \completeer{\Cantorspace[n]}
      \forall s \in \Cantorspace[n^a - (n + 1)] \ \psi^a_n(r, s)
        \extendedby \psi^\gamma_n(r, s)$.
  \end{enumerate}
  An approximation $a$ is \definedterm{$X'$-terminal} if no
  configuration is compatible with both $X'$ and a one-step extension
  of $a$. Let $\calA(a, X')$ denote the family of sets of the form
  $\saturation{\image{(\phi_X \composition \phi^\gamma)}
  {\Cantorspace[n^a]}}{E \intersection F}$, where $\gamma$ varies
  over all configurations compatible with $a$ and $X'$.
  
  \begin{lemma} \label{dichotomy:terminalimpliesintersecting}
    Suppose that $X' \subseteq X$ and $a$ is an $X'$-terminal
    approximation. Then $\calA(a, X')$ is $E$-locally $F$-intersecting.
  \end{lemma}
  
  \begin{lemmaproof}
    Suppose, towards a contradiction, that there are configurations
    $\gamma_0$ and $\gamma_1$, both compatible with $a$ and $X'$,
    such that $\image{(\phi_X \composition \phi^{\gamma_0})}
    {\Cantorspace[n^a]}$ and $\image{(\phi_X \composition \phi^
    {\gamma_1})}{\Cantorspace[n^a]}$ are $F$-disjoint sets contained in
    the same $E$-class. Fix $f \from \completeer{\Cantorspace[n^a]} \to
    \Bairespace$ such that $(\phi_{E \setminus F} \composition
    f)(r) = \sequence{(\phi_X \composition \phi^{\gamma_i})(r(i))}[i < 2]$
    for all $r \in \completeer{\Cantorspace[n^a]}$, and let
    $\gamma$ be the configuration given by $n^\gamma = n^a + 1$,
    $\phi^\gamma(s \concatenation \sequence{i}) = \phi^{\gamma_i}(s)$
    for all $i < 2$ and $s \in \Cantorspace[n^a]$, $\psi_n^\gamma(r, s
    \concatenation \sequence{i}) = \psi_n^{\gamma_i}(r, s)$ for all $i <
    2$, $n < n^a$, $r \in \completeer{\Cantorspace[n]}$, and $s \in
    \Cantorspace[n^a - (n + 1)]$, and $\psi_{n^a}^\gamma(r,
    \emptysequence) = f(r)$ for all $r \in \completeer{\Cantorspace
    [n^a]}$. Then $\gamma$ is compatible with a one-step extension of
    $a$, contradicting the fact that $a$ is $X'$-terminal.
  \end{lemmaproof}
  
  Proposition \ref{dichotomy:intersecting} and Lemma \ref
  {dichotomy:terminalimpliesintersecting} ensure that if $a$ is
  $B^\alpha$-terminal, then there is an $(E \intersection F)$-invariant
  \Borel partial quasi-transversal $B(a, B^\alpha)$ of $E$ over $E
  \intersection F$ puncturing $\calA(a, B^\alpha)$. Let $B^{\alpha + 1}$
  be the set obtained from $B^\alpha$ by subtracting the union of the
  sets of the form $B(a, B^\alpha)$, where $a$ varies over all
  $B^\alpha$-terminal approximations.
  
  \begin{lemma}
    \label{dichotomy:nonterminalapproximationextension}
    Suppose that $\alpha < \omega_1$ and $a$ is a non-$B^{\alpha +
    1}$-terminal approximation. Then $a$ has a
    non-$B^\alpha$-terminal one-step extension.
  \end{lemma}
  
  \begin{lemmaproof}
    Fix a one-step extension $b$ of $a$ for which there is a
    configuration $\gamma$ compatible with $b$ and $B^{\alpha + 1}$.
    Then $\image{(\phi_X \composition \phi^\gamma)}{\Cantorspace
    [n^b]} \subseteq B^{\alpha + 1}$, so $b$ is not $B^\alpha$-terminal. 
  \end{lemmaproof}
  
  Fix $\alpha < \omega_1$ such that the families of $B^\alpha$- and
  $B^{\alpha + 1}$-terminal approximations coincide, and let $a_0$
  denote the unique approximation for which $n^{a_0} = 0$. As $\calA
  (a_0, X') = \set{\equivalenceclass{x}{E \intersection F}}[x \in X']$ for all
  $X' \subseteq X$, we can assume that $a_0$ is not $B^\alpha$-terminal,
  since otherwise $B^{\alpha+1} = \emptyset$, so $X$ is a countable union
  of $(E \intersection F)$-invariant \Borel partial quasi-transversals of $E$
  over $E \intersection F$.
  
  By recursively applying Lemma \ref
  {dichotomy:nonterminalapproximationextension}, we obtain
  non-$B^\alpha$-terminal one-step extensions $a_{n+1}$ of $a_n$
  for all $n \in \N$. Define $\phi \from \Cantorspace \to \Bairespace$ by
  $\phi(c) = \union[n \in \N][\phi^{a_n}(\restriction{c}{n})]$, as well as
  $\psi_n \from \completeer{\Cantorspace[n]} \times \Cantorspace \to
  \Bairespace$ by $\psi_n(r, c) = \union[m > n][\psi_n^{a_m}(r,
  \restriction{c}{(m - (n + 1))})]$ for all $n \in \N$. Clearly these functions are
  continuous.
  
  \begin{lemma}
    The function $\phi_X \composition \phi$ is a homomorphism from
    $\Ezero \setminus \diagonal{\Cantorspace}$ to $E \setminus F$.
  \end{lemma}
  
  \begin{lemmaproof}
    We will show that if $c \in \Cantorspace$, $n \in \N$, and $r \in
    \completeer{\Cantorspace[n]}$, then
    \begin{equation*}
      (\phi_{E \setminus F} \composition \psi_n)(r, c) = \sequence{(\phi_X
        \composition \phi)(r(i) \concatenation \sequence{i} \concatenation
          c)}[i < 2].
  \end{equation*}
  As $X \times X$ is \Hausdorff, it is sufficient to show that if $U$ is an
  open neighborhood of $\sequence{(\phi_X \composition \phi)(r(i)
  \concatenation \sequence{i} \concatenation c)}[i < 2]$ and $V$ is an
  open neighborhood of $(\phi_{E \setminus F} \composition \psi_n)(r,
  c)$, then $U \intersection V \neq \emptyset$. Towards this end, fix
  $m > n$ such that $\product[i < 2][\image{\phi_X}{\extensions{\phi^
  {a_m}(r(i) \concatenation \sequence{i} \concatenation s)}}] \subseteq
  U$ and $\image{\phi_{E \setminus F}}{\extensions{\psi_n^{a_m}(r,
  s)}} \subseteq V$, where $s = \restriction{c}{(m - (n + 1))}$. As
  $a_m$ is not $B^\alpha$-terminal, there is a configuration $\gamma$
  compatible with it. Then $\sequence{(\phi_X \composition \phi^\gamma)
  (r(i) \concatenation \sequence{i} \concatenation s)}[i < 2] \in U$ and
  $(\phi_{E \setminus F} \composition \psi_n^\gamma)(r, s) \in V$, thus $U
  \intersection V \neq \emptyset$.
  \end{lemmaproof}
  
  Set $\phi' = \phi_X \composition \phi$. As the equivalence relations $E' =
  \preimage{(\phi' \times \phi')}{E}$ and $F' = \preimage{(\phi' \times \phi')}
  {F}$ have the \Baire property, Proposition \ref{generalization:meager}
  ensures that they are meager, in which case the closed equivalence
  relation $D' = \preimage{(\phi' \times \phi')}{\diagonal{X}}$ is nowhere
  dense, so Proposition \ref{generalization:homomorphism} yields a
  continuous homomorphism $\psi \from \Cantorspace \to \Cantorspace$
  from $\triple{\setcomplement{\diagonal{\Cantorspace}}}{\Ezero}
  {\setcomplement{\Ezero}}$ to $\triple{\setcomplement{D'}}{\Ezero}
  {\setcomplement{(E' \union F')}}$. As $\phi'$ is a homomorphism from
  $\Gzero$ to $E$, it follows that $\Gzero \subseteq E'$, so Proposition \ref
  {generalization:generates} ensures that $\Ezero \subseteq E'$, thus $\phi'
  \composition \psi$ is a continuous embedding of $\pair{\Ezero}{\diagonal
  {\Cantorspace}}$ into $\pair{E}{F}$.
\end{theoremproof}

\begin{remark} \label{dichotomy:DC}
  The apparent use of choice beyond $\DC$ in the above argument
  can be eliminated by first running the analog of the argument using
  the weakening of Proposition \ref{dichotomy:intersecting} without
  any definability constraints on the partial quasi-transversal puncturing
  the family (which can be proven in the same manner, but without
  using Proposition \ref{generalization:separation}), in order to obtain
  an upper bound $\alpha' < \omega_1$ on the least ordinal $\alpha <
  \omega_1$ for which the sets of $B^\alpha$- and
  $B^{\alpha+1}$-terminal approximations coincide.
\end{remark}

A \definedterm{partial quasi-uniformization} of a set $R \subseteq X \times
Y$ over an equivalence relation $F$ on $Y$ is a set $S \subseteq R$ for
which there exists $k \in \N$ such that every vertical section of $S$ is
contained in a union of at most $k$ $F$-classes. Theorem \ref
{intro:uniformization} is a consequence of Proposition \ref
{generalization:transversal:two}, Theorems \ref
{generalization:uniformization} and \ref{finitebasis:uniformization}\footnote
{Theorem \ref{finitebasis:uniformization} can be replaced with the usual
\Lusin--\Novikov uniformization theorem to establish the special case
referred to in Remark \ref{intro:uniformization:F}.}, and the following:

\begin{theorem} \label{dichotomy:quasiuniformization}
  Suppose that $X$ and $Y$ are \Hausdorff spaces, $E$ is an analytic
  equivalence relation on $X$, $F$ is a \Borel equivalence relation on
  $Y$, and $R \subseteq X \times Y$ is an $(E \times \diagonal
  {Y})$-invariant analytic set whose vertical sections are contained in
  countable unions of $F$-classes. Then exactly one of the following
  holds:
  \begin{enumerate}
    \item There is a cover $\sequence{R_n}[n \in \N]$ of $R$ by $(\restriction
      {(E \times F)}{R})$-invariant \Borel-in-$R$ partial quasi-uniformizations
      of $R$ over $F$.
    \item There are continuous embeddings $\pi_X \from \Cantorspace \into
      X$ of $\Ezero$ into $E$ and $\pi_Y \from \Cantorspace \into Y$ of
      $\diagonal{\Cantorspace}$ into $F$ such that $\image{(\pi_X \times
      \pi_Y)}{\Ezero} \subseteq R$.
  \end{enumerate}
\end{theorem}

\begin{theoremproof}
  To see that the conditions are mutually exclusive, note that if both hold,
  then there exists $n \in \N$ for which $\image{\projection[1]}{\preimage
  {(\pi_X \times \pi_Y)}{R_n}}$ is a non-meager partial quasi-transversal of
  $\Ezero$ with the \Baire property, a contradiction.

  To see that at least one of the conditions holds, observe that $\restriction
  {(E \times \completeer{Y})}{R}$ is analytic and $\restriction{(\completeer
  {X} \times F)}{R}$ is \Borel, and appeal to Theorem \ref
  {dichotomy:quasitransversal} to see that if condition (1) fails, then there is
  a continuous embedding $\pi \from \Cantorspace \into R$ of $\pair{\Ezero}
  {\diagonal{\Cantorspace}}$ into $\pair{\restriction{(E \times \completeer
  {Y})}{R}}{\restriction{(E \times F)}{R}}$. It follows that the function $\pi_X'
  = \projection[X] \composition \pi$ is a continuous reduction of $\Ezero$ to
  $E$ and the function $\pi_Y' = \projection[Y] \composition \pi$ is a
  continuous homomorphism from $\Ezero \setminus \diagonal
  {\Cantorspace}$ to $\setcomplement{F}$, and therefore from $\Gzero$ to
  $\setcomplement{F}$, so Proposition \ref{generalization:meager} ensures
  that the relation $F' = \preimage{(\pi_Y' \times \pi_Y')}{F}$ is meager, in
  which case the closed relations $D_X' = \preimage{(\pi_X' \times \pi_X')}
  {\diagonal{X}}$ and $D_Y' = \preimage{(\pi_Y' \times \pi_Y')}{\diagonal
  {Y}}$ are nowhere dense, so Proposition \ref
  {generalization:homomorphism} yields a continuous homomorphism $\pi'
  \from \Cantorspace \to \Cantorspace$ from $\triple{\setcomplement
  {\diagonal{\Cantorspace}}}{\Ezero}{\setcomplement{\Ezero}}$ to $\triple
  {\setcomplement{(D_X' \union D_Y')}}{\Ezero}{\setcomplement{(\Ezero
  \union F')}}$, thus the functions $\pi_X = \pi_X' \composition \pi'$ and
  $\pi_Y = \pi_Y' \composition \pi'$ are continuous embeddings of $\Ezero$
  into $E$ and $\diagonal{\Cantorspace}$ into $F$. As $\image{\pi}
  {\Cantorspace} \subseteq R$, it follows that $\image{(\pi_X' \times \pi_Y')}
  {\diagonal{\Cantorspace}} \subseteq R$, so the facts that $\pi_X'$ is a
  homomorphism from $\Ezero$ to $E$ and $R$ is $(E \times \diagonal
  {Y})$-invariant ensure that $\image{(\pi_X' \times \pi_Y')}{\Ezero}
  \subseteq R$, thus the fact that $\pi'$ is a homomorphism from $\Ezero$
  to $\Ezero$ implies that  $\image{(\pi_X \times \pi_Y)}{\Ezero} \subseteq
  R$.
\end{theoremproof}

We say that a set $R \subseteq X \times X$ is the \definedterm{graph of a
partial quasi-function} over an equivalence relation $F$ on $X$ if there
exists $k \in \N$ such that $\forall \sequence{x_i, y_i}[i < k] \in R^k \ (\forall
i, j < k \ x_i \mathrel{F} x_j \implies \exists i < j < k \ y_i \mathrel{F} y_j)$.
We say that $R$ is the \definedterm{graph of a partial quasi-injection} if
both $R$ and $\inverse{R}$ have this property. Theorem \ref{intro:graph}
is a consequence of Proposition \ref{generalization:bijection}, Theorem
\ref{finitebasis:graph}, and the following:

\begin{theorem} \label{dichotomy:quasigraph}
  Suppose that $X$ is a \Hausdorff space, $E$ is an analytic equivalence
  relation on $X$, $F$ is a \Borel equivalence relation on $X$, and every
  $E$-class is a countable union of $(E \intersection F)$-classes. Then
  exactly one of the following holds:
  \begin{enumerate}
    \item There is a cover $\sequence{R_n}[n \in \N]$ of $E \setminus F$ by
      $((E \intersection F) \times (E \intersection F))$-invariant \Borel-in-$E$
      graphs of partial quasi-injections over $E \intersection F$.
    \item There is a continuous embedding $\pi \from \Cantorspace \times 2
      \into X$ of $\pair{\Ezero \times \completeer{2}}{\Ezero \disjointunion
      \diagonal{\Cantorspace}}$ into $\pair{E}{F}$.
  \end{enumerate}
\end{theorem}

\begin{theoremproof}
  To see that the conditions are mutually exclusive, note that if both hold,
  then there exists $n \in \N$ with the property that $\set{c \in \Cantorspace}
  [\pi(c,0) \mathrel{R_n} \pi(c,1)] \times \set{1}$ is a non-meager partial
  quasi-transversal of $\Ezero \times \diagonal{\set{1}}$ with the \Baire
  property, a contradiction.

  To see that at least one of the conditions holds, note first that if there are
  $\restriction{((E \intersection F) \times F)}{(E \setminus F)}$-invariant
  \Borel-in-$E$ partial quasi-uniformizations $R_n$ of $E$ over $E
  \intersection F$ for which $E = \union[n \in \N][R_n]$, then the sets of the
  form $R_m \intersection \inverse{R_n}$, where $m, n \in \N$, are as
  desired. By Theorem \ref{dichotomy:quasiuniformization}, we can
  therefore assume that there are continuous embeddings $\pi_X \from
  \Cantorspace \into X$ of $\Ezero$ into $E$ and $\pi_Y \from \Cantorspace
  \into X$ of $\diagonal{\Cantorspace}$ into $F$ such that $\image{(\pi_X
  \times \pi_Y)}{\Ezero} \subseteq E \setminus F$. As $\Ezero$ has
  countable index below the equivalence relation $E' = \preimage{(\pi_X
  \times \pi_X)}{E}$, the latter is meager, so the closed subequivalence
  relation $D' = \preimage{(\pi_X \times \pi_X)}{\diagonal{X}}$ is nowhere
  dense, thus Proposition \ref{generalization:homomorphism} yields a
  continuous homomorphism $\pi' \from \Cantorspace \to \Cantorspace$
  from $\triple{\setcomplement{\diagonal{\Cantorspace}}}{\Ezero}
  {\setcomplement{\Ezero}}$ to $\triple{\setcomplement{D'}}{\Ezero}
  {\setcomplement{E'}}$. Define $\pi \from \Cantorspace \times 2 \to X$ by
  $\pi(c, 0) = (\pi_X \composition \pi')(c)$ and $\pi(c, 1) = (\pi_Y \composition
  \pi')(c)$ for all $c \in \Cantorspace$.
\end{theoremproof}

\begin{acknowledgements}
  We would like to thank Alexander Kechris for asking the questions that
  led to this work, as well as Julia Millhouse for bringing several typos to
  our attention.
\end{acknowledgements}

\bibliographystyle{amsalpha}
\bibliography{bibliography}

\end{document}